\documentclass[11pt,a4paper,amsart]{article}
\usepackage{amsmath,amsfonts,amsbsy,latexsym,epsfig,amsthm,amssymb,graphics,eepic,graphicx}

\input xy
\xyoption{all}

\title{Hurwitz Spaces and Moduli Spaces as Ball Quotients via Pull-back}

\author{Brent R. Doran}

\date{}

\newtheorem{Lemma}{Lemma}
\newtheorem{Thm}{Theorem}
\newtheorem{Prop}{Proposition}
\newtheorem{Cor}{Corollary}
\newtheorem{Def}{Definition}

\theoremstyle{remark}

\theoremstyle{remark}
\newtheorem{Rmk}{Remark}

\newcommand{\pone}{\ensuremath{\mathbb{P}^1}}
\newcommand{\PmS}{\ensuremath{\mathbb{P}^1 \setminus S}}

\newcommand{\Ldual}{L^\vee}

\newcommand{\DM}{\ensuremath{DM(n,\mu)}}

\newcommand{\Spi}{\ensuremath{\mathcal{S}_{\pi}}}

\else\oddsidemargin25mm\evensidemargin25mm\marginparwidth25mm\fi
\begin{document}

\maketitle

\abstract{  We define hypergeometric functions using intersection
homology valued in a local system.  Topology is emphasized;
analysis enters only once, via the Hodge decomposition.  By a
pull-back procedure we construct special subsets $\Spi$, derived
from Hurwitz spaces, of Deligne-Mostow moduli spaces $\DM$.
Certain $\DM$ are known to be ball quotients, uniformized by
hypergeometric functions valued in a complex ball (i.e., complex
hyperbolic space).  We give sufficient conditions for $\Spi$ to be
a subball quotient. Analyzing the simplest examples in detail, we
describe ball quotient structures attached to some moduli spaces
of inhomogeneous binary forms. This recovers in particular the
structure on the moduli space of rational elliptic surfaces given
by Heckman and Looijenga. We make use of a natural partial
ordering on the Deligne-Mostow examples (which gives an easy way
to see that the original list of Mostow, eventually corrected by
Thurston, is in error), and so highlight two key examples, which
we call the Gaussian and Eisenstein ancestral examples.}

\section{Introduction}\label{sec:Intro}

    A number of classical moduli spaces $\mathcal{M}$ admit the structure of a
locally symmetric space $\Gamma \backslash G/K$, where $K$ is a
maximal compact subgroup and $\Gamma$ is a discrete subgroup of
$G$.  The identification is given by a $\Gamma$-invariant map
$\Phi: G/K \rightarrow \mathcal{M}$ that descends to an
isomorphism from $\Gamma \backslash G/K$ to $\mathcal{M}$.

\vspace{.1in}

\xymatrix{ & & & & &  & & {G/K} \ar[dl]_{\Phi} \ar[d]^{/ \sim} \\
& & & & & &  {\mathcal{M} } & {\ar[l]_{\cong} \Gamma \backslash
G/K} }

\vspace{.1in}

In the language of the 19th century, $(G/K, \Phi)$ is a {\em
uniformization} of $\mathcal{M}$. If $\Gamma$ does not act freely,
it is an {\em orbifold uniformization}. For example, when $G =
PU(1, n)$ then $G/K$ is the complex $n$-ball $\mathbb{B}^n$, or
complex hyperbolic $n$-space, and so the uniformization endows
$\mathcal{M}$ with a complex hyperbolic metric (possibly with
orbifold singularities). We call $\mathcal{M} \cong \Gamma
\backslash \mathbb{B}^n$ a ``ball quotient".

In their seminal work on hypergeometric functions, Deligne and
Mostow \cite{DM,Mos1,DM2} proved that certain geometric invariant
theory (GIT) moduli spaces $\DM$ of $n$ points on $\pone$ are ball
quotients. Here the uniformizing group $\Gamma$ is a monodromy
representation of the spherical braid group on $n$ strings, the
ball is $\mathbb{B}^{n-3}$, and the uniformizing map $\Phi$ is the
single-valued inverse to a map $HG_{\mu}$ built from multi-valued
generalized hypergeometric functions on $\DM$.

\vspace{.1in}

\xymatrix{ & & & & &   & \ \ \mathbb{B}^{n-3} \ar[d]^{\ / \sim} \ar@{->}@(l,u)[1,-1]_{\Phi_{\mu}} \\
& & & & & {\DM} \ar@3{->}[ur]^{HG_{\mu}} \ar[r]^{\ \ \cong} &
{\Gamma \backslash \mathbb{B}^{n-3}} }

\vspace{.1in}

In this paper, we make use of an overlooked, essentially
topological, property of hypergeometric functions defined via
local systems. This yields a method for producing subball
quotients of $\DM$ that have natural geometric interpretations in
terms of moduli.

Intersection homology has several convenient attributes
\cite{Borel,GM1,GM2,MacP}. One is that it comes with an
intersection pairing whose signature can be computed using
explicit cycles. Here, this pairing defines a Hermitian form,
$\Psi$, which only depends on a list $\mu$ of $n$ fractions
(associated with the $n$ points on $\pone$). Consequently, varying
the coordinates of the points, i.e., moving through the moduli
space, preserves the form.  In particular, the monodromy group
$\Gamma$ preserves $\Psi$: so the data $\mu$ determine a {\em
lattice} $\Lambda$ over a ring of integers $R$, together with
$\Gamma$ acting as automorphisms of $\Lambda$. Thus $\mu$ in fact
defines a Hermitian locally symmetric space. This procedure
parallels the main approach, using ordinary and compactly
supported cohomologies, of Deligne and Mostow.

Also, intersection homology is a bivariant functor that it is
insensitive to points of trivial local monodromy. Consequently,
given $\pi: \pone \rightarrow \pone$, and a rank $1$ local system
$l_T \rightarrow \pone \setminus T$, the intersection homology of
the pull-back rank 1 local system $\pi^*(l_T) \rightarrow \pone
\setminus \pi^{-1}(T)$ only ``sees'' the interesting points $S
\subset \pi^{-1}(T)$ that contribute to monodromy.  This is the
feature that allows us to relate certain Hurwitz spaces to subball
quotients.

Finally, intersection homology here admits a Hodge decomposition
into orthogonal subspaces corresponding to holomorphic and
anti-holomorphic forms. This decomposition encodes the only actual
analysis to enter into what is otherwise a topological
construction. We define the multi-valued map $HG_{\mu}$ to be the
coordinate expression of the holomophic section of a flat
Grassmannian bundle (over the moduli space of $n$ distinct points
on $\pone$) given by the linear subspace of holomorphic 1-forms.
For the Deligne-Mostow examples, the flat Grassmanian bundle is
actually a flat projective space bundle and $\Psi$ is of
Lorentzian signature, so that $HG_{\mu}$ is valued in a complex
ball in projective space, as desired.

Section \ref{sec:HGfunctions} develops the theory of
hypergeometric functions of Deligne-Mostow type from the
perspective of intersection homology valued in a local system.  We
take a hands-on approach to understanding the cycles, the form
$\Psi$, and the monodromy group $\Gamma$.  No claim is made as to
original results, although some proofs may be new. The reader is
referred to the original paper of Deligne and Mostow \cite{DM} for
some technical details that carry over to intersection homology
{\em mutatis mutandis}, and for a complete discussion of the
theory of holomorphic and anti-holomorphic forms valued in a rank
1 local system on an $n$-punctured $\pone$. The section concludes
by recalling the fundamental uniformization results of Deligne and
Mostow \cite{DM,Mos1}, via conditions $INT$ and $\Sigma INT$.

Section \ref{Sec:Ancestral} serves two purposes, both related to
the $\Psi$-lattices that underlie the Deligne-Mostow
uniformizations. First it establishes that two important examples,
which we call the Eisenstein and Gaussian ancestral examples, have
a rich supply of subball quotients. This is key for the results of
Section \ref{sec:Pullback}. Second, the section makes concrete a
fact which is implicit in Deligne and Mostow's analysis, but is
often overlooked. Namely, GIT-stable collisions of points on
$\pone$ in a Deligne-Mostow uniformization correspond to orbifold
subball quotients (basically, the mirrors of complex reflections
associated to the monodromy group $\Gamma$).  This gives an easy
way to organize much of the Deligne-Mostow list of
uniformizations; for instance, when $n \geq 7$ every example but
one is a collision ``descendant'' of an ancestral example. In
particular, it is immediate that Mostow's original list
\cite{Mos1} is missing a number of examples. (Thurston corrected
that list by computer in \cite{Thur}, and to the best of our
knowledge his list is complete, although a proof is not provided.)
We show in \cite{Dor1} (and it is also seen by somewhat different
means in \cite{DKvG}) that one of these ``missing'' examples is a
cover of the ball quotient structure on the moduli space of cubic
surfaces from \cite{ACT}.

Section \ref{sec:Pullback} discusses the pull-back procedure.
Given a map $\pi: \pone \rightarrow \pone$ and a collection of
points $T$, let $\nu$ denote the data of the non-trivial local
monodromy of a rank $1$ local system $l_T \rightarrow \pone
\setminus T$.  Thus $\nu$ may be written as a list of $|T|$
fractions.  The pull-back local system $\pi^*l_T \rightarrow \pone
\setminus \pi^{-1}(T)$ is specified by its local monodromy list of
fractions $\mu$.  Given $\nu$, the data $\mu$ is determined
entirely by the ramification of $\pi$ over $T$.  Consequently,
varying $\pi$ while preserving the ramification over $T$
(describing a Hurwitz space of $\pone$ covers of $\pone$), results
in a constrained variation of the points of $\pi^{-1}(T)$. The
variation of the subset $S \subset \pi^{-1}(T)$ of non-trivial
local monodromy points thus determines a subvariety $\Spi$ of a
Deligne-Mostow moduli space $\DM$. If this moduli space is a ball
quotient, then one may ask if $\Spi$ is a subball quotient via
restriction of hypergeometric functions. We give a sufficiency
criterion. Detailed classification in very special cases with $|T|
= 3$ produces a list of $\Spi$ that at once are codimension 1
subball quotients and admit natural finite covers by moduli spaces
of inhomogeneous binary forms.  The codimension 1 subball quotient
of the Eisenstein ancestral example may be interpreted as the
moduli space of rational elliptic surfaces, via the Weierstrass
fibration description of Miranda \cite{Mir}, thus recovering the
ball quotient structure on that space described in \cite{HL}.

Throughout we assume the number of points on $\pone$ is $n \geq
3$.

\section{Hypergeometric Functions after Deligne and
Mostow}\label{sec:HGfunctions}

\subsection{Background on local systems and subsystems}\label{subsec:LocalSys}

Let $X$ be a connected manifold. The following are equivalent
characterizations of a complex local system $L \rightarrow X$ up
to isomorphism:
\begin{enumerate}
\item complex vector bundle with flat connection

\item locally constant sheaf of complex vector spaces

\item $\pi_1(X,0)$ representation on a complex vector space, known
as the {\em monodromy} representation, where $0$ is some chosen
base point in $X$
\end{enumerate}

To pass from the first description to the second, identify the
flat vector bundle with its sheaf of locally constant sections.
Fix a base point $0 \in X$.  Then the functor ``fiber at $0$"
produces a complex vector space endowed with an action of
$\pi_1(X,0)$, hence deriving the third description from the
second.  Note that the choice of base point does not affect the
isomorphism class, since a different base point yields the same
monodromy representation up to conjugacy in the general linear
group of the vector space.

In the case of rank $1$ local systems, the monodromy
representation is one-dimensional hence abelian.  Since $H_1(X)$
is naturally isomorphic to $\pi_1(X)^{ab}$, the local system is
determined by a homomorphism $H_1(X) \rightarrow \mathbb{C}^*$.

The simplest case is also the basic object of study for this
paper, namely, rank $1$ local systems on the projective line
punctured at $n$ points $s_j, j \in \{1, \ldots, n \}$.  Observe
$H_1(\PmS)$ is an abelian group.  Take positively oriented circles
centered at the $s_j$ to be representative cycles of the group
generators.  The only relation is that the product of the
generators is the identity.  It is easy to show that:

\begin{Prop}\label{Prop:rk1localmonod}
Given a set of points $S = \{ s_1, \ldots, s_n \}$ in $\pone$ and
a set of complex numbers ${\mu} = \{ \mu_1, \ldots, \mu_n \}$,
there is, up to isomorphism, a unique rank $1$ local system on
$\PmS$ with $(S, e^{2 \pi i {\mu}})$ as the local monodromy data.
However, the local system is {\em not} determined up to {\em
unique} isomorphism; the fibers may be uniformly rescaled by any
element of $\mathbb{C}^*$.
\end{Prop}

If the $\mu_j$ are real then $\alpha_j = e^{2 \pi i \mu_j}$ is on
the unit circle in $\mathbb{C}^*$.  In particular the local
monodromy may be of finite order.  For the purposes of the
Deligne-Mostow theory of hypergeometric functions, one takes
$\mu_j \in \mathbb{Q}, \forall j$.  Ultimately this condition will
follow from the constraints $INT$ or $\Sigma INT$ on $\mu$ that
guarantee uniformization, so it is a matter of convenience to
demand it in advance.
\begin{Def} \label{Def:DM-Lsystem}
A {\em Deligne-Mostow local system} is a rank $1$ local system $L
\rightarrow \PmS$, where $S$ is a finite collection of points,
such that for all $j, \mu_j$ (encoding the local monodromy data
$\alpha_j = e^{2 \pi i \mu_j}$ at $s_j$) is a rational number.
Clearly the $\mu_j$ may be adjusted to lie between $0$ and $1$
without changing $L$.
\end{Def}

By a {\em local subsystem} of a local system $L$ we mean a locally
constant subsheaf.  Note that the monodromy data of a DM local
system is defined over the ring of integers $R =
\mathbb{Z}[\zeta_d]$, where $\zeta_d$ is a $d^{th}$ root of unity,
and $d$ is the least common denominator of the $\mu_j$.  We call
the corresponding local subsystem with fiber $R$ the {\em
Deligne-Mostow local subsystem}, denoted $L(R)$.

The {\em dual local system} $\Ldual$ will be needed for the
homology theory of the next section.  It has a straightforward
explicit description.

\begin{Prop}\label{Prop:Ldual}
If $L$ is the DM local system determined by the data $(S,
{\alpha})$, equivalently by $(S, {\mu})$, then $\Ldual$ is the DM
local system determined by $(S, \overline{{\alpha}})$,
equivalently by $(S, 1-{\mu})$.  In other words, $\Ldual =
\overline{L}$.  Furthermore $\Ldual(R) = \overline{L}(R)$, where
$R$ is the ring of integers defining the Deligne-Mostow local
subsystem of $L$.
\end{Prop}
\begin{proof}
By Proposition \ref{Prop:rk1localmonod}, both $L$ and $\Ldual$ are
characterized by local monodromy data.  Thus, if the data $(S, \{
\alpha_j \})$ determine $L$, then $\Ldual$ is characterized by
$(S, \{\alpha_j^{-1} \})$. Furthermore, if $L$ is a Deligne-Mostow
local system then $\alpha_j^{-1} = \overline{\alpha}_j$, because
$\alpha_j$ lies on the unit circle in $\mathbb{C}$.  So $\Ldual$
is determined by $\overline{\alpha}_j$, and hence by $-\mu_j$, or
equivalently, by $1-\mu_j$ (normalizing to lie between $0$ and
$1$).  It is immediate from the explicit local monodromy data that
all of the corresponding DM local subsystems are defined over the
same ring of integers $R$.
\end{proof}

\subsection{Intersection homology valued in Deligne-Mostow local
systems}

\subsubsection{Background on intersection homology}

What follows is an informal discussion.  The goal is to impart
intuition and to highlight the results needed in the sequel.
Details for the trivial local system case can be found in
\cite{Borel,GM1,GM2}, and the arguments are easily adapted for
general local systems (see also \cite{MacP}).

Intersection homology can be defined for any Whitney stratified
pseudo-manifold.  Any quasi-projective variety $X$ admits a
Whitney stratification, where the unique open stratum,
$X^{nonsing}$ is the ``nonsingular part" of $X$. Intersection
homology is a topological invariant, independent of the choice of
stratification.  The simplest definition is the original
formulation, due to Goresky and MacPherson, in terms of geometric
chains.  Many models for the chains are acceptable, but for our
purposes piecewise linear chains are perfectly satisfactory.

Intersection homology theory is similar to ordinary homology
theory on $X$.  The boundary operator is the same but the
intersection chain complex, $IC_{\cdot}(X)$ is a subcomplex of the
ordinary chain complex.  Those ordinary chains whose intersection
with the singular locus $X^{sing}$ are too ``perverse", i.e., too
non-generic, are disallowed.  A choice $\overline{p}$ of {\em
perversity} is then a choice of which chains are admissible. The
default choice of perversity for algebraic varieties is ``middle
perversity."  Middle perversity intersection homology of $X$,
denoted $IH_*(X)$, has many nice properties --- the so-called {\em
K\"{a}hler package}.  The most important property here is
Poincar\'{e}-Verdier duality:

\begin{Prop} \label{Prop:P-V-duality} (Poincar\'{e}-Verdier duality)
Let $d$ be the real dimension of $X$. There is a non-degenerate
bilinear pairing
$$IH_k(X) \otimes IH_{d-k}(X) \stackrel{\cap}{\longrightarrow}
\mathbb{C}$$
\end{Prop}

Let us now be more precise.  Let $X^{sing}$ be stratified by
$\{S_{\beta}\}$, where $\beta$ is the codimension of $S_{\beta}$
in $X$.

\begin{Def}\label{Def:perversity}
\begin{itemize}
\item A (classical) {\em perversity} $\overline{p}$ is a positive
integer-valued non-decreasing function on the natural numbers
$\{2, \ldots, dim_{\mathbb{R}}(X) \}$, satisfying $\overline{p}(2)
= 0$ and $\overline{p}(\beta + 1) \leq \overline{p}(\beta) + 1$.

\item An $i$-chain $\xi$ in $X$ is an {\em intersection $i$-chain}
if it satisfies the {\em admissibility} conditions:
    \begin{enumerate}
        \item $dim_{\mathbb{R}}( \xi \cap S_{\beta}) \leq
        dim_{\mathbb{R}}(S_{\beta}) + dim_{\mathbb{R}}(\xi) - n +
        \overline{p}(\beta)$
        \item $dim_{\mathbb{R}}( \partial \xi \cap S_{\beta}) \leq
        dim_{\mathbb{R}}(S_{\beta}) + dim_{\mathbb{R}}(\partial \xi) - n +
        \overline{p}(\beta)$
    \end{enumerate}
\end{itemize}
\end{Def}

\begin{Rmk}
The perversity starts with codimension $2$ because the singular
locus of a pseudo-manifold is real codimension at least $2$. The
second admissibility condition ensures the intersection chains
form a complex.
\end{Rmk}

Just as with homology, intersection homology can be valued in
sheaves other than the constant sheaf. In particular one considers
intersection homology valued in a complex local system $L
\rightarrow X \setminus S$.  The standard notation (using middle
perversity) is $IH_*(X, L)$.  Here $X^{sing}$ is $S$.  The support
of an $L$-valued intersection chain is, as before, a geometric
chain in $X$ that satisfies admissibility conditions based on
choice of perversity. The only difference is that, in the
nonsingular locus $X^{nonsing}$, the ``value" attached to the
chain is a section of $L$ over the chain. In other words, an
$L$-valued intersection chain is an ordinary chain in
$X^{nonsing}$, valued in $L$, for which the closure of its support
satisfies the admissibility criteria to be an intersection chain
in $X$.

When $L = {\mathbb{C}}$, the trivial rank $1$ local system on $X$,
one of course recovers the usual intersection homology with
complex coefficients.

\subsubsection{Vector space structure with basis}

Now let $L \rightarrow \PmS$ be a DM local system.

\begin{Lemma} \label{Lem:badboundary}
The geometric support of $IC_0(\pone, L)$ and of $IC_1(\pone, L)$
is $\pone \setminus S$.  Let $K \subset S$ be the subset of points
of nontrivial monodromy.  Then the geometric support of
$IC_2(\pone, L)$ is $\pone \setminus K$.
\end{Lemma}
\begin{proof}
All choices of (classical) perversity are equivalent for a
one-dimensional complex variety, because $\overline{p}(2) = 0$. So
the middle perversity is the zero perversity $\overline{p} \equiv
0$. Consequently, the first admissibility criterion disallows both
points $s_j \in S$ and any $1$-chains that intersect an $s_j$, but
imposes no constraint on the $2$-chains. The second admissibility
criterion does not restrict the $0$- and $1$-chains further.

The application of the second admissibility criterion to
$2$-chains is more subtle, because the chains are valued in a
local system.  There are two types of intersection with $S$: the
$2$-chain either contains an $s_j$ with non-trivial monodromy
($\alpha_j \neq 1$) or it only contains $s_j$ with trivial
monodromy.  In the latter case, the boundary of the $2$-chain is a
$1$-cycle that encloses but does not intersect $s_j$, and so is
admissible.  In the former case, observe that any such $L$-valued
2-chain has as boundary a $1$-cycle that intersects $s_j$ and so
is {\em not} admissible.
\end{proof}

\vspace{.1in}

\noindent {\bf Example:} Let $C$ denote a small circle oriented
counter-clockwise and centered at $s_k$, and let $\theta$ denote
the line segment from $s_k$ to $p$ on $C$.  Denote the choice of
section of $L$ at $p$ by $\hat{p}$, and its horizontal extension
over $C$ and $\theta$ by $\hat{C}$ and $\hat{\theta}$,
respectively.  This determines a unique horizontal section over
the disk $D$, denoted $\hat{D}$, with discontinuities (when
$\alpha_k \neq 1$) along $\theta$.

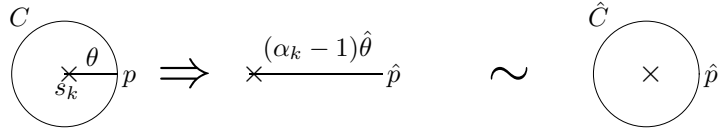
\vspace{.1in}

\begin{picture}(320,50)

\put(60,25){\circle{40}}  \put(57,22){\large $\times$}

\put(39,43){$C$}  \put(56,17){$s_k$}

\put(60,25){\line(1,0){20}}  \put(68,28){$\theta$}

\put(82,22){$p$}  \put(95,20){\huge $\Rightarrow$}


\put(130,25){\line(1,0){50}}  \put(127,22){\large $\times$}
\put(182,22){$\hat{p}$}  \put(135,32){$(\alpha_k - 1)
\hat{\theta}$}

\put(220,20){\huge $\sim$}

\put(280,25){\circle{40}}  \put(258,43){$\hat{C}$}
\put(277,22){\large $\times$}  \put(302,22){$\hat{p}$}

\end{picture}

\vspace{.1in}

\noindent It is easy to see that $\partial \hat{D} = \hat{C} -
(\alpha_k -1) \hat{\theta}$.  In particular, the support of the
boundary intersects $s_k$, violating the second admissibility
criterion.

Intersection homology is insensitive to points $s_j$ of trivial
local monodromy. More precisely, any intersection $1$-cycle
enclosing such a point $s_j$ is homologous to an intersection
$1$-cycle that does not enclose it.  This homology is realized by
an intersection $2$-chain that contains $s_j$ and takes values in
the trivial local system (extended over $s_j$). Formalizing this
argument yields the following Lemma.

\begin{Lemma}\label{Lem:IH-insensitive-trivialmonod}
Let $L \rightarrow \PmS$ be a rank $1$ complex local system.  Let
$K$ be the subset of points $\{ s_{m_1}, \ldots, s_{m_k} \}
\subset S$ with nontrivial local monodromy, that is, those points
$s_j$ with $\alpha_j \neq 1$. Let $\hat{L}$ denote the local
system on $\pone \setminus K$ defined by the local monodromies
$\alpha_{m_j}$. Then there is a natural isomorphism $IH_1(\pone,
L) \cong IH_1(\pone, \hat{L})$.
\end{Lemma}

Thus the study of intersection homology valued in Deligne-Mostow
local systems reduces to considering those local systems defined
by $\alpha_j \neq 1, \forall j$. In that case, it is elementary to
prove that the first homology groups in all the usual homology
theories are isomorphic.

\begin{Lemma}\label{Lemma:allhomologiesequal}
Let $K$ and $\hat{L}$ be as above.  Then there are natural
isomorphisms
$$IH_1(\pone, \hat{L}) \cong H_1(\pone \setminus K, \hat{L}) \cong
H_1^{lf}(\pone \setminus K, \hat{L}),$$ where $H_1^{lf}$ denotes
locally finite homology.
\end{Lemma}

\begin{Prop}\label{Prop:IHbetti}
Let $K = \{s_{m_1}, \ldots, s_{m_k}\}$ be the subset of points
$s_j$ in $S$ with $\alpha_j \neq 1$, and define $\hat{L}
\rightarrow \pone \setminus K$ as before.  Assume $K \neq
\emptyset$. Then
$$ dim_{\mathbb{C}}(IH_0(\pone, \hat{L})) = 0 \, , \
dim_{\mathbb{C}}(IH_1(\pone, \hat{L})) = k-2 \, , \
dim_{\mathbb{C}}(IH_2(\pone, \hat{L})) = 0 \ . $$
\end{Prop}

\begin{proof}

 Use the
isomorphisms of Lemma \ref{Lem:IH-insensitive-trivialmonod} and
Lemma \ref{Lemma:allhomologiesequal} to identify $IH_1(\pone, L)$
with $H_1^{lf}(\pone \setminus K, \hat{L})$. A good choice of $k$
generators for $H_1^{lf}(\pone \setminus K, \hat{L})$ is the set
of $1$-cycles $\gamma_j$ with endpoints $s_{m_j}, s_{m_{j+1}}$
(where the final one, $\gamma_k$, connects endpoints $s_{m_k},
s_1$ in that order). Without loss of generality, let the $s_{m_j}$
be aligned along the equator, so that the $\gamma_j$ themselves
form the equator.

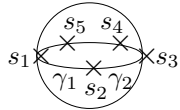
\vspace{.1in}

\begin{picture}(320,50)

\put(170,25){\circle{40}}

\put(170,25){\ellipse{40}{10}}

\put(147,22){\large $\times$}  \put(187,22){\large $\times$}
\put(167,17.5){\large $\times$}

\put(157,27){\large $\times$}  \put(177,27){\large $\times$}
\put(160,35){$s_5$}  \put(174,35){$s_4$}

\put(139,22){$s_1$}  \put(168,10){$s_2$}  \put(195,22){$s_3$}

\put(156,14){$\gamma_1$}  \put(177,14){$\gamma_2$}


\end{picture}

\vspace{.1in}


Relations among the homology generators are precisely those linear
combinations of $\gamma_j$ which are the boundary of some locally
finite $2$-chain $\hat{D}$.  It is clear that the support of
$\hat{D}$ must be either the upper or lower hemisphere.  Pick a
point $p$ in the upper hemisphere (the choice of hemisphere is not
important) and a section $\hat{p}$ of $\hat{L}$ over $p$. The
section has a unique horizontal extension over the hemisphere
containing $p$, so that $\partial \hat{D} = \sum_j \gamma_j$, so
in homology
$$\sum_j \gamma_j = 0.$$
This is the first linear relation among the generators. Picking a
different lift than $\hat{p}$ simply rescales the section, and
hence the boundary relation, by a complex number; so this choice
doesn't alter the linear relation. The horizontal section extends
to the lower hemisphere, but now a choice must be made: the
natural extension is to a {\em multi-section} on $\pone \setminus
K$.  The choice therefore lies in selecting a $\gamma_j$ over
which to {\em continuously} extend the horizontal section into the
lower hemisphere, to get a single-valued section; but it is not
continuous along the remaining $\gamma_i, i \neq j$.  By crossing
at $\gamma_j$, the resulting $\hat{D}$ (with $D$ now the lower
hemisphere) has a boundary that can be explicitly written in terms
of local monodromies, yielding the second linear relation on the
homology generators:
$$
\alpha_2^{-1} \alpha_3^{-1} \cdots \alpha_j^{-1} \gamma_1 + \ldots
\alpha_j^{-1} \gamma_{j-1} + \gamma_j + \alpha_{j+1} \gamma_{j+1}
+ \ldots \alpha_{j+1} \cdots \alpha_k \gamma_k = 0.
$$
Observe that a different choice of $\gamma_j$ simply rescales the
linear combination by a complex number, so there is no change to
the relation.  Thus there are precisely two relations on the $k$
generators, so $dim_{\mathbb{C}}(IH_1(\pone, \hat{L})) = k-2$.

The other intersection homologies are easy to compute.
$IH_0(\pone, \hat{L}) = 0$ because for any point $p \in \pone
\setminus K$ and section $\hat{p}$, the boundary of a 1-cycle
whose support passes through $p$ that loops around precisely one
$s_{m_j}$ is just $(\alpha_{m_j} - 1)\hat{p}$, and so $\hat{p}$ is
a boundary of an intersection $1$-cycle. $IH_2(\pone, \hat{L}) =
0$ because any $2$-chain has non-trivial boundary so there are no
$2$-cycles.
\end{proof}

\begin{Cor} \label{Cor:IH-grouping-into-two-sets}
Given any two points $s_i, s_j \in K$, the intersection $1$-chain
$I_{i,j} \equiv \frac{1}{\alpha_i - 1}\hat{C}_i + \hat{\gamma} +
\frac{1}{1-\alpha_j}\hat{C}_j$ is in fact a $1$-cycle not
homologous to zero.  Moreover, given any partition of a subset of
$K$ into two disjoint collections $\{s_i \}_{i \in I}$ and $\{s_j
\}_{j \in J}$ $(I, J$ index sets such that $I \cap J = \emptyset)$
where $\prod_{i \in I} \alpha_i \neq 1$ and $\prod_{j \in J}
\alpha_j \neq 1$ (or equivalently, $\sum_{i \in I} \mu_i$ and
$\sum_{j \in J} \mu_j \not\in \mathbb{Z}$), the analogous
$1$-chain $I_{I,J}$ that encircles the two collections and
connects them by a segment $\gamma$ is in fact a $1$-cycle.
\end{Cor}

\begin{proof}
The locally finite $1$-cycle with support a line segment
$\gamma_{i,j}$ from $s_i$ to $s_j$ and section determined by
extending $\hat{p}$ at $p$ is in the same locally finite homology
class as the intersection $1$-cycle $I_{i,j} \equiv
\frac{1}{\alpha_i - 1}\hat{C}_i + \hat{\gamma} +
\frac{1}{1-\alpha_j}\hat{C}_j$. In particular, it is non-zero in
locally finite homology.  Because $IC_2(\pone, \hat{L}) \subset
C_2^{lf}(\pone \setminus K, \hat{L})$, it is clear $I_{i,j}$
cannot be the boundary of any intersection $2$-chain, and so is
non-zero in intersection homology.

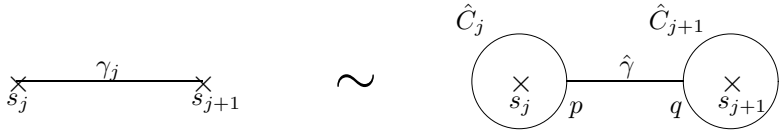
\vspace{.2in}

\begin{picture}(320,50)


\put(16,21){\large $ \times $}  \put(16,16){$s_j$}

\put(20,25){\line(1,0){70}}

\put(50,28){$\gamma_j$}

\put(86,21){\large $\times $}  \put(86,16){$s_{j+1}$}

\put(140,20){\huge $\sim$}


\put(210,25){\circle{36}}  \put(206,21){\large $\times$}
\put(206,15){$s_j$}

\put(290,25){\circle{36}}  \put(286,21){\large $\times$}
\put(285,15){$s_{j+1}$}

\put(228,25){\line(1,0){44}}


\put(229,13){$p$}  \put(267,13){$q$}

\put(248,28){$\hat{\gamma}$}

\put(186,45){$\hat{C}_j$}  \put(259,45){$\hat{C}_{j+1}$}

\end{picture}

\vspace{.2in}

The statement for collections of points $\{ s_i \}_{i \in I}$ and
$\{ s_j \}_{j \in J}$ is immediate. Denote the closed curve that
encircles the first collection (and no other $s_l$) by $C_1$, mark
a point $p \in C_1$ and a choice of section at $p$ by $\hat{p}$
that extends to a section $\hat{C}_1$. Likewise about the second
collection construct $C_2$ and $q$, and connect $p$ and $q$ with
the line segment $\gamma$.  Extend the horizontal section from
$\hat{p}$ to $\hat{\gamma}$ and $\hat{C}_2$.

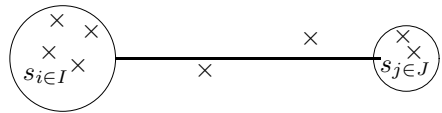
\vspace{.1in}
\begin{picture}(320,80)

\put(80,40){\circle{40}}

\put(87,48){$\times$}  \put(71,40){$\times$} \put(82,35){$\times$}
\put(74,52){$\times$}  \put(65,32){$s_{i \in I}$}

\put(100,40){\line(1,0){100}}  \put(130,33){$\times$}
\put(170,45){$\times$}

\put(210,40){\circle{25}}

\put(200,35){$s_{j \in J}$}  \put(205,46){$\times$}
\put(209,40){$\times$}

\end{picture}
\vspace{.1in}

The boundary of the intersection $1$-chain
$$ I_{I,J} =
\frac{1}{\prod_{i \in I} \alpha_i - 1}\hat{C}_1 + \hat{\gamma} +
\frac{1}{1-\prod_{j \in J} \alpha_j}\hat{C}_2
$$
is zero, hence it is a $1$-cycle.
\end{proof}

We therefore get an intersection homology basis taken from the set
$\{I_{i,i+1} \}_{i \in \{1, \ldots, k-1\}} \cup \{ I_{k,1} \}$.

\begin{Cor}\label{Cor:IHbasis}
Partition $S$ into two subsets $S_1$ and $S_2$, whose elements are
indexed by $i$ and $j$ respectively.  The cycles $\{I_{i,i+1}
\equiv \frac{1}{\alpha_i - 1}\hat{C}_i + \hat{\gamma} +
\frac{1}{1-\alpha_{i+1}}\hat{C}_{i+1} \}_{i \in \{1, \ldots,
|S_1|-1 \}}$ taken together with the cycles $\{ I_{j,j+1} \}_{j
\in \{|S_1|+1, \ldots, |S_2|-1 \}}$ form a basis for $IH_1(\pone,
L)$ if and only if $\sum_{i \in S_1} \mu_i \not\in \mathbb{Z}$ (or
equivalently $\sum_{j \in S_2} \mu_j \not\in \mathbb{Z}$).
\end{Cor}

\begin{Rmk}
In fact, more is true.  Given any partition of $S$ into $S_1$ and
$S_2$ (satisfying the above condition on the $\mu_i$ and $\mu_j$),
{\em any} tree connecting the points of $S_1$ taken together with
any tree connecting the points of $S_2$ defines a basis in locally
finite homology.  This translates into a basis for intersection
homology in the fashion indicated above. For a proof of the
locally finite homology fact, see \cite[Section 2.5]{DM}
\end{Rmk}

\subsubsection{Intersection pairing: skew-Hermitian form on
$IH_1(\pone,L)$}\label{subsec:Psi}

\begin{Prop}\label{Prop:(S,mu)-determine-Psi}
The data $(S, \mu)$ determine, up to a real scalar, the
intersection pairing on intersection homology.  The pairing puts a
skew-Hermitian form on $IH_1(\pone, L)$, so multiplication by the
complex number $\imath = \sqrt{-1}$ yields an Hermitian form
$\Psi$, unique up to a real scalar.
\end{Prop}

\begin{proof}
Poincar\'{e}-Verdier duality gives a nondegenerate bilinear
pairing between $IH_k(X,L)$ and $IH_{d-k}(X, \Ldual)$, where $d$
is the real dimension of $X$ and $\Ldual$ is the dual local
system.  By Proposition \ref{Prop:Ldual}, for a Deligne-Mostow
local system $\Ldual = \overline{L}$. Because $IH_*(\overline{L})
= \overline{IH_*(L)}$, the duality pairing is
$$IH_k(X, L) \otimes IH_{d-k}(X, \Ldual) = IH_k(X, L) \otimes
\overline{IH_{d-k}(X,L)} \stackrel{\cap}{\longrightarrow} L
\otimes \Ldual \cong {\mathbb{C}} $$ Thus there is a
skew-Hermitian intersection form on $IH_1(X,L)$.  From this, a
Hermitian form $\Psi$ is obtained by multiplication by $\imath$.
By Proposition \ref{Prop:rk1localmonod}, given $\{\mu_j\}$, $L$ is
determined up to a $\mathbb{C}^*$ factor, so the intersection
pairing and Hermitian form $\Psi$ are determined up to a real
scalar.
\end{proof}

The intersection pairing on intersection $1$-cycles, expressed in
the basis of Corollary \ref{Cor:IHbasis}, is a skew-Hermitian
matrix $Int$. Let $Int(i,j)$ denote the $(i,j)$ entry of $Int$. If
$|i-j| > 1$ then $Int(i,j) = 0$, because the support of
$I_{i,i+1}$ doesn't intersect that of $I_{j,j+1}$. It remains to
compute the self-intersection of $I_{i, i+1}$ and the intersection
number for adjacent basis cycles (when $|i-j| = 1$).

\begin{Prop}\label{Prop:IntMatrix}
The skew-Hermitian intersection form, with respect to the basis
$I_{i,i+1}, i \in \{1, \ldots, k-2 \}$, is the matrix $Int$ with
entries:
$$
Int(i,j) = \left\{ \begin{array}{ccl}
                    \frac{1}{1- \alpha_i} - 1 +
                    \frac{1}{1- \alpha_{i+1}} & , & j = i   \\
                    -\frac{1}{1-\alpha_i} & , & j = i+1  \\
                     \frac{1}{1-\overline{\alpha}_i} & , & j=i-1 \\
                     0 & , &        |j - i| > 1
                 \end{array}  \right.
$$

\end{Prop}

\begin{proof}

The computation is immediate from the following picture.  The
positive orientation is taken to be counterclockwise.

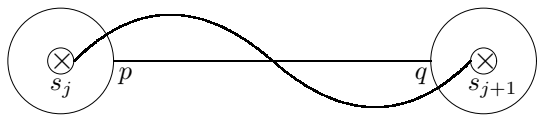
\vspace{.1in}
\begin{picture}(320,50)

\put(80,25){\circle{40}}  \put(240,25){\circle{40}}

\put(76,22){\large $\times$}  \put(236,22){\large $\times$}

\put(100,25){\line(1,0){120}}

\put(102,18){$p$}   \put(214,18){$q$}

\qbezier(85,25)(120,60)(160,25)

\qbezier(160,25)(200,-10)(235,25)

\put(80,25){\circle{10}}   \put(240,25){\circle{10}}

\put(76,14){$s_j$}   \put(234,14){$s_{j+1}$}

\end{picture}

\vspace{.1in}

Note that the deformation chosen to compute the intersection
number is particularly convenient given the choice of section
(with discontinuities at $p$ and $q$). A different choice would,
of course, yield the same number, albeit presented as a sum of
different terms.

\end{proof}

The Hermitian form $\Psi$ is simply $\imath$ times the
intersection pairing, so in matrix form, $\Psi(j,k) = \imath
Int(j,k)$.

\vspace{.1in}

\noindent {\bf Example:}
 When $n = 4$,
$$
\Psi = \imath \left[  \begin{array}{cc}
                \frac{1}{1- \alpha_1} - 1 +
                    \frac{1}{\alpha_2}  &
                    -\frac{1}{1-\alpha_2} \\
                \frac{1}{1-\overline{\alpha}_2} &
                \frac{1}{1- \alpha_2} - 1 +
                    \frac{1}{\alpha_3}
                \end{array}
        \right]
$$
In particular, if $\mu_i = \frac{1}{2}, \forall i$, then
$$
\Psi = \left[
            \begin{array}{cc}
            0           & -\frac{\imath}{2} \\
            \frac{\imath}{2} &   0
            \end{array}
        \right]
$$

\vspace{.1in}

One application is that the signature of the form can be computed
purely in terms of $\sum_i \mu_i$.  There are a number of ways to
show this.  We give a constructive argument, which produces an
explicit basis for a maximal positive definite subspace and its
orthogonal negative definite subspace in $IH_1(\pone, \hat{L})$.

Let $I_{I,J}$ be the intersection $1$-cycle described above,
enclosing $s_i, i \in I$ with $C_1$ and $s_j, j \in J$ with $C_2$.
Recall that the $\mu_i \in (0,1)$.  The monodromy along $C_1$ is
given by $\prod_{i \in I} \alpha_i$, and so is determined by the
fractional part of $\sum_{i \in I} \mu_i$.  The analogous
statement holds for the monodromy along $C_2$.

\begin{Lemma} \label{Lem:Psi-Signature-prelims}
Let $Frac(x)$ denote the fractional part of the non-negative real
number $x$, i.e., $x-\lfloor x \rfloor$. If the sum $Frac(\sum_{i
\in I} \mu_i) + Frac(\sum_{j \in J} \mu_j) < 1$ then the length of
$I_{I,J}$ is negative, that is, $\Psi(I_{I,J}, I_{I,J}) < 0$. If
$Frac(\sum_{i \in I} \mu_i) + Frac(\sum_{j \in J} \mu_j) = 1$ then
$\Psi(I_{I,J}, I_{I,J}) = 0$. If $Frac(\sum_{i \in I} \mu_i) +
Frac(\sum_{j \in J} \mu_j) > 1$ then $\Psi(I_{I,J}, I_{I,J}) > 0$.
\end{Lemma}
\begin{proof}
This is just clever work with the self-intersection number
computed above:
$$ -\frac{1}{\prod_{i \in I} \alpha_i - 1} -1 + \frac{1}{1 - \prod_{j \in
J}\alpha_j}$$
\end{proof}

\begin{Prop}\label{Prop:signature}
The signature of $\Psi$ is $(\sum \mu_i -1, \sum (1 - \mu_i) -1)$
where the first term is the dimension of a maximal positive
definite subspace and the second is the dimension of a maximal
negative definite subspace.
\end{Prop}

\begin{proof}
The strategy is to build a succession of $I_{I,J}$ which are
mutually orthogonal.  Since the sign of the $\Psi$-length of
$I_{I,J}$ is known by Lemma \ref{Lem:Psi-Signature-prelims}, we
get explicit maximal positive definite and negative definite
subspaces.  We inductively construct $I_{I,\{j\}}$ as shown in the
following picture.

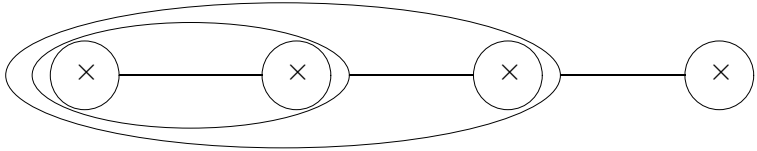
\vspace{.1in}

\begin{picture}(320,60)

\multiput(40,30)(80,0){4}{\circle{26}}

\put(53,30){\line(1,0){54}}

\multiput(36 ,28)(80,0){4}{\large $\times$}


\put(80,30){\ellipse{120}{40}}

\put(140,30){\line(1,0){47}}


\put(115,30){\ellipse{210}{55}}

\put(220,30){\line(1,0){47}}

\end{picture}

\vspace{.1in}

In particular, assuming none of these is zero length, for each
positive integer less than $\sum_i \mu_i$ we produce a new
positive length vector in the positive definite subspace, linearly
independent from the preceding $I_{I,J}$.  All of the remaining
basis vectors generated by this procedure are in the orthogonal
negative definite subspace. Because the total dimension is $n-2$,
we get the stated result.

Now assume that some $I_{I,\{j\}}$ so constructed has zero length.
Select $a$ and $b$ so that $aI_{I,\{j\}} + bI_{\{j\},\{j+1\}}$ has
positive length.  Then $\overline{a}I_{I,\{j\}} + \overline{b}
I_{\{j,j+1\}}$ is orthogonal and has negative length.  Furthermore
these are orthogonal to all previously constructed vectors in the
sequence.  In particular the number of positive and negative
definite vectors produced remains unchanged.  Continue the
inductive procedure as before.

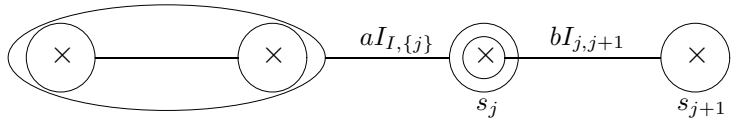
\vspace{.1in}

\begin{picture}(320,60)

\multiput(40,30)(80,0){4}{\circle{26}}

\put(53,30){\line(1,0){54}}

\multiput(36 ,28)(80,0){4}{\large $\times$}


\put(140,30){\line(1,0){47}}


\put(200,30){\circle{16}}

\put(153,35){$aI_{I,\{j\}}$}   \put(225,35){$b I_{j,j+1}$}

\put(197,10){$s_j$}  \put(273,10){$s_{j+1}$}


\put(208,30){\line(1,0){59}}

\put(80,30){\ellipse{120}{40}}

\end{picture}

\vspace{.1in}

\end{proof}

\subsubsection{Lattice structure over ring of integers}

We now recall the notion of a lattice from the theory of modules.
The ring $R$ will always denote a ring of algebraic integers.  We
will study the structure in greater depth in Section
\ref{Sec:Ancestral}.

\begin{Def}\label{Def:Lattice}
A module-theoretic {\em lattice} $\Lambda$ is a finite rank module
over a ring $R$, endowed with an Hermitian form $\Psi$ taking
values in $R$.  A lattice is {\em unimodular} if the determinant
of $\Psi$ (the {\em discriminant} of $\Lambda$) is $\pm 1$.  A
{\em sublattice} $M$ of $\Lambda$ is a submodule together with the
restriction of $\Psi$. A sublattice $M$ is {\em primitive} if
there is no other sublattice $M^{\prime}$ such that $M = r
M^{\prime}$ for $r$ not a unit in $R$. The group $Aut(\Lambda)$ of
{\em unitary transformations} (equivalently, lattice
automorphisms) of $\Lambda$ is the group of module automorphisms
of $\Lambda$ that preserve $\Psi$.
\end{Def}

The intersection homology valued in a DM local system has the
structure of a lattice.

\begin{Lemma}\label{Lem:ringofint}
Let $R$ be the ring of integers in $\mathbb{Q}[\zeta_d]$, where
$d$ is the lowest common denominator of the $\mu_i$, and $\zeta_d$
is a primitive $d^{th}$ root of unity. Then $\Psi$ is defined over
$R$.
\end{Lemma}
\begin{proof}
It is immediate from our matrix descriptions of $\Psi$ that it is
defined over $R$.  More formally, this follows by considering the
local subsystem $L(R)$ with fiber the subring $R \subset
\mathbb{C}$. The pairing of $L(R)$ with $L^{-1}(R) =
\overline{L(R)}$ induces an $R$-valued skew-Hermitian pairing on
$IH_1(\pone, L(R))$.  This pairing may be identified with the
pairing from Section \ref{subsec:Psi}, since it was there only
determined up to a real scalar by the data $\mu$.
\end{proof}

\subsection{Configuration space of $n$ points on $\mathbb{P}^1$}

Now let the positions of the $n$ points $\{ s_j \}$ vary on
$\pone$, while fixing the $\{ \mu_j \}$.  We start with some basic
definitions.

\begin{Def}\label{Def:Config-Braid}
\begin{itemize}
\item Let $\mathcal{P}_n$ denote the {\em configuration space} of
$n$ distinct ordered points on $\pone$.  That is, $\mathcal{P}_n =
(\pone)^n \setminus \{z_i = z_j, i \neq j \}$.

\item Let $\mathcal{P}^{\Sigma}_n := \mathcal{P}_n/\Sigma_n$,
where $\Sigma_n$ is the symmetric group on $n$ letters which acts
by permuting the $s_j$. This is the configuration space of $n$
unordered points on $\pone$.  It is an open subset of
$\mathbb{P}^n$.

\item We refer to $\pi_1(\mathcal{P}_n)$ as the {\em braid group}
on $n$ colored strings on $\pone$.
\end{itemize}

\item Similarly, $\pi_1(\mathcal{P}_n/\Sigma_n)$ is the braid
group on $n$ colorless strings on $\pone$.
\end{Def}

The automorphism group of $\pone$ is $PGL_2(\mathbb{C})$.  An
automorphism is completely determined by its action on any three
distinct points of $\pone$. The diagonal action on $(\pone)^n$
restricts to a free action on $\mathcal{P}_n$.

\begin{Def}\label{Def:moduli}
Let $\mathcal{M}_n$ denote the {\em moduli space} of $n$ distinct
ordered points on $\pone$.  It is the quotient of $\mathcal{P}_n$
under the free diagonal action of $PGL_2(\mathbb{C})$.  That is,
$\mathcal{M}_n \equiv \mathcal{P}_n/PGL_2(\mathbb{C}) \cong (\pone
\setminus \{0,1,\infty\})^{n-3} \setminus \{z_i = 0,1,\infty, z_j,
i \neq j \}$.  Similarly, the moduli space of $n$ distinct
unordered points on $\pone$ is the $\Sigma_n$ quotient of
$\mathcal{M}_n$, which we denote by $\mathcal{M}^{\Sigma}_n$.
\end{Def}

Observe that $\mathcal{P}_n \cong \mathcal{M}_n \times
PGL_2(\mathbb{C})$. Consequently, ignoring the choice of base
point because the spaces in question are all connected, we have
\begin {Lemma}\label{Lem:relate-Mn-and-Pn}
\begin{eqnarray*}\pi_1(\mathcal{P}_n) = \pi_1(\mathcal{M}_n) \times
\pi_1(PGL_2(\mathbb{C})) = \pi_1(\mathcal{M}_n) \times
\mathbb{Z}/2\mathbb{Z} \\ \pi_1(\mathcal{P}^{\Sigma}_n) =
\pi_1(\mathcal{M}^{\Sigma}_n) \times \pi_1(PGL_2(\mathbb{C})) =
\pi_1(\mathcal{M}^{\Sigma}_n) \times \mathbb{Z}/2\mathbb{Z}
\end{eqnarray*}
\end{Lemma}

As discussed in Section \ref{subsec:LocalSys}, local systems on
$\mathcal{P}_n$ are characterized up to isomorphism by
representations of $\pi_1(\mathcal{P}_n, x_0)$ (the choice of
$x_0$ is irrelevant up to isomorphism).  The braid group acts
naturally on the intersection homology of $L \rightarrow \PmS$,
and thus defines a local system $\mathcal{L}$ of rank $k-2$ on
$\mathcal{P}_n$.  (See Section \ref{Subsec:Monodromy} for a
detailed discussion of the braid representation.)

A point $p \in \mathcal{P}_n$ specifies a subset $S(p) \subset
\pone$ of $n$ distinct points.  Let $L_p \rightarrow \pone
\setminus S(p)$ denote the DM local system defined by the data
$(S(p), \mu)$.  It turns out that what one would hope for is in
fact true: namely, the vector spaces $IH_1(\pone, L_p)$ arrange
themselves over $\mathcal{P}_n$ into the local system $\mathcal{L}
\rightarrow \mathcal{P}_n$.  There is some ambiguity because $L_p$
is not determined up to unique isomorphism by $\mu$ (Proposition
\ref{Prop:rk1localmonod}). More precisely, one finds (adapting the
arguments from cohomology to intersection cohomology is immediate
when $\alpha_i \neq 1, \forall i$) \cite[pp. 22, 26-27]{DM}:

\begin{Lemma}\label{Cor:Lexists-almostunique}
Given $\mu$ there is a local system $\mathcal{L} \rightarrow
\mathcal{P}_n$ with fiber at $p$ given by $H^1(\pone \setminus
S(p), L_p)$, where $L_p$ is the DM local system with monodromy
data $\mu$.  This local system is unique up to tensor product with
a rank one local system $\mathcal{O} \rightarrow \mathcal{P}_n$.
\end{Lemma}

The ambiguity is removed by projectivizing the fibers. One may
think of the resulting {\em canonical} flat projective bundle
$P\mathcal{L}$ as one of a number of canonical Grassmanian
bundles, constructed in the analogous way, on $\mathcal{M}_n$.
More precisely, let $dim_{\mathbb{C}}(IH_1(\pone, L_p)) = k-2$.
For each $l, 1 < l < k-2$, there is a canonical flat Grassmannian
bundle $Gr(l, \mathcal{L})$, with fiber $Gr(l,k-2)$ over each $m
\in \mathcal{M}_n$.

Such bundles are characterized by their monodromy representation.
\begin{Def}
Denote the monodromy group of $\mathcal{L}$ by $\Gamma$, and the
(canonical given $\mu$) monodromy group associated to
$P\mathcal{L}$ by $P\Gamma$.
\end{Def}

Furthermore, because the $PGL_2(\mathbb{C})$ action is trivial on
the projective fibers, this flat bundle of projective spaces
descends to $\mathcal{M}_n$.  Alternatively, the projective
representation of $\pi_1(\mathcal{P}_n)$ encoded in the
projectivization of the local system $\mathcal{L}$ is
simultaneously a projective representation for
$\pi_1(\mathcal{M}_n)$ by Lemma \ref{Lem:relate-Mn-and-Pn}, and so
canonically describes a flat projective space bundle on
$\mathcal{M}_n$.  To be more precise, let $\Theta$ denote the
projective monodromy group of $P\mathcal{L} \rightarrow
\mathcal{M}_n$.  We verify that $P\Gamma = \Theta$.
\begin{Prop}\label{Prop:PGamma}
The projective monodromy group $\Theta$ of the flat bundle of
projective spaces $P\mathcal{L} \rightarrow \mathcal{M}_n$ is
isomorphic to the projective monodromy group $P\Gamma$ of
$P\mathcal{L} \rightarrow \mathcal{P}_n$.
\end{Prop}
\begin{proof}
The flat bundle of projective spaces $P\mathcal{L} \rightarrow
\mathcal{M}_n$ is isomorphic to a flat subbundle of the bundle of
projective spaces $\mathbb{P}(IH_1(\pone, L_p)) \rightarrow
\mathcal{P}_n$: simply restrict the bundle via the inclusion $i:
\mathcal{M}_n \subset \mathcal{P}_n, m = (m_0, \ldots, m_{n-3})
\mapsto (0,1,\infty, m_0, \ldots, m_3)$.  Thus the projective
monodromy representation $\Theta(i_*\pi_1(\mathcal{M}_n))$ of
$P\mathcal{L}$ is automatically a subgroup of $P\Gamma =
\Theta(\pi_1(\mathcal{P}_n))$. By Lemma
\ref{Lem:relate-Mn-and-Pn}, $\pi_1(\mathcal{P}_n) =
\pi_1(\mathcal{M}_n) \times \mathbb{Z}/2\mathbb{Z}$ so
$\Theta(\pi_1(\mathcal{P}_n))$ is isomorphic to
$\Theta(i_*(\pi_1(\mathcal{M}_n)))$ twisted by a character of
$\mathbb{Z}/2\mathbb{Z}$.  In particular, they define equivalent
projective representations.
\end{proof}
\begin{Rmk}
This justifies using $P\Gamma$ in either context, so henceforth we
will not refer to $\Theta$, only to $P\Gamma$.  We freely use
whichever interpretation is convenient, without further comment,
throughout.
\end{Rmk}

\subsection{Definition of hypergeometric functions}

In the preceding sections, all of the results were topological.
Analysis enters via the definition of hypergeometric functions.

\begin{Prop}\label{Prop:IH-Hodge}
Let $L \rightarrow \PmS$ be a DM local system.  There is an
orthogonal ``Hodge decomposition" $IH_1(\pone, L) \cong
IH_{1,0}(\pone, L) \oplus IH_{0,1}(\pone, L)$.  The Hermitian form
$\Psi$ on $IH_1(\pone, L)$ is positive definite on the subspace
$IH_{1,0}(\pone, L)$ and negative definite on $IH_{0,1}(\pone,
L)$.
\end{Prop}
\begin{proof}
In general the decomposition follows from work of Saito
\cite{Saito1,Saito2} when $L$ is a local system of geometric
origin in the sense of Grothendieck-Deligne (i.e., is a polarized
variation of Hodge structure). It can be seen more directly by
interpreting $IH_{1,0}$ as the space of holomorphic $L$-valued
1-forms and $IH_{0,1}$ as the space of holomorphic
$\overline{L}$-valued 1-forms (i.e., anti-holomorphic $L$-valued
1-forms). See \cite[Section 2]{DM} for details: the isomorphism of
homology theories when $\alpha_s \neq 1$, from Lemma
\ref{Lemma:allhomologiesequal}, together with the non-degeneracy
of $\Psi$, carry over their argument unchanged.
\end{proof}

\begin{Cor}
\begin{eqnarray*}
dim_{\mathbb{C}}(IH_{1,0}(\pone,L)) = (\sum_i \mu_i) - 1 \\
dim_{\mathbb{C}}(IH_{0,1}(\pone,L)) = (\sum_i 1 - \mu_i) -1
\end{eqnarray*}
\end{Cor}

\begin{Cor}
If $\sum_i \mu_i = 2$ then $\Psi$ is signature $(1,n-3)$.
\end{Cor}

Fix $\mu$ such that $\sum_i \mu_i = 2$.  Pick coordinates on the
fiber of $P\mathcal{L} \rightarrow \mathcal{M}_n$ at some point
$m_0$ and extend by flatness.  We use the fact that
$IH_{1,0}(\pone, L)$ is a distinguished 1-dimensional subspace,
spanned by some $\omega_{\mu}$ to define:
\begin{Def}
The multi-valued holomorphic map $HG_{\mu}: \mathcal{M}_n
\rightarrow \mathbb{P}^{k-2}$ is the coordinate expression of
$\omega_{\mu}$.   We call it the {\em $\mu$-hypergeometric
function of Deligne-Mostow type}.
\end{Def}
By construction, $HG_{\mu}(m)$ is an orbit of the projective
monodromy group of $P\mathcal{L} \rightarrow \mathcal{M}_n$.
$HG_{\mu}$ is completely determined up to automorphisms of
$\mathbb{P}^{k-2}$.

\begin{Rmk}
Let $\Sigma$ denote the symmetries of the list $\mu = (\mu_0,
\ldots, \mu_{n-1})$.  Then $\Sigma$ acts  on $\mathcal{M}_n$ as
permutations of the associated coordinates $s_i$.  It is clear
from the definition that $HG_{\mu}$ descends to a map from
$\mathcal{M}^{\Sigma}_n$. We denote this map by $HG_{\mu}$ as
well.  The domain will always be clear from context.
\end{Rmk}

\begin{Rmk}
This definition of hypergeometric functions may be generalized to
arbitrary $\mu$ by using the coordinate expression for
$IH_{1,0}(\pone,L)$ in the corresponding flat Grassmannian bundle.
I am not aware of an analog of this definition in the literature.
\end{Rmk}

\begin{Rmk}
When $\sum_i \mu_i = 2$ as above, there is in fact a unique
holomorphic 1-form up to scaling.  It may be written as:
$$
\omega_{\mu} = \prod_i (z-s_i)^{-\mu_i} e \cdot dz,
$$
where $e$ is a horizontal multi-section of $L$ (to cancel the
monodromy of the function so that $\omega_{\mu}$ is a well-defined
section).  This is the famous hypergeometric 1-form.
\end{Rmk}

\subsection{Uniformization by a complex ball}

\subsubsection{Complex ball and discrete subgroups of $PU(1,n)$}

Let $\Psi$ be an Hermitian form of signature $(1,n)$ on
$\mathbb{C}^{n+1}$.

\begin{Def}\label{Def:Ball}
The {\em complex ball} $\mathbb{B}^n \subset \mathbb{P}^n$ is
defined to be the subset of points that lift to vectors in
$\mathbb{C}^{n+1}$ of strictly positive $\Psi$-length.  In
particular, $\Psi$ defines a {\em complex hyperbolic} metric on
$\mathbb{B}^n$.
\end{Def}

\begin{Rmk}
An Hermitian form over $\mathbb{C}$ is determined up to
equivalence (change of coordinates) by its signature.
Consequently, $\mathbb{B}^n$ is independent of $\Psi$.  In
particular, $\mathbb{B}^n$ is a ``ball" because, in an appropriate
coordinate system $z = (z_0, \ldots, z_n)$,
\begin{eqnarray*}
\Psi(z, z) = |z_0|^2 - |z_1|^2 - \cdots |z_n|^2 > 0 \\
\Rightarrow 1 > |\frac{z_1}{z_0}|^2 + |\frac{z_2}{z_0}|^2 + \cdots
|\frac{z_n}{z_0}|^2
\end{eqnarray*}
\end{Rmk}

Let $PU(1,n)$ denote the group of projective linear
transformations that lift to linear transformations on
$\mathbb{C}^{n+1}$ which preserve $\Psi$.

$$
\begin{array}{cccl}
U(1,n) & \hookrightarrow & GL_{n+1}(\mathbb{C}) &  =
Aut(\mathbb{C}^{n+1}) \\
\downarrow & & \downarrow & \\
PU(1,n) &  \hookrightarrow & PGL_{n+1}(\mathbb{C}) & =
Aut(\mathbb{P}^n)
\end{array}
$$

The complex ball has an interpretation as an Hermitian symmetric
domain of type $I$.
\begin{Prop}\label{Prop:BallAsSymmSp}
The complex ball is a symmetric space.
$$\mathbb{B}^{n} \cong PU(1,n)/P(U(1) \times U(n-1))$$
Furthermore $Aut(\mathbb{B}^n) \cong PU(1,n)$.
\end{Prop}
\begin{proof}
See \cite{KN}, Volume II, Example 10.7, pages 282--285.
\end{proof}

\begin{Prop}\label{Prop:PU(1,n)}
When $\sum \mu_i = 2$, the projective monodromy group $P\Gamma$ of
$P\mathcal{L}$ is a subgroup of $PU(1,n-3)$ and so acts as
automorphisms of the complex ball $\mathbb{B}^{n-3}$.
\end{Prop}
\begin{proof}
The braid group acts via compactly supported isotopies on $\pone$.
These isotopies induce automorphisms of the local system $L_{p_0}
\rightarrow \pone \setminus S(p_0)$.  Any such automorphism is
multiplication in the fibers by $\mathbb{C}^*$ (by Proposition
\ref{Prop:rk1localmonod}). This in turn induces a constant
$\mathbb{C}^*$ rescaling of $IH_1(\pone, L_{p_0})$, and so a real
rescaling of the skew-Hermitian intersection pairing.  Hence up to
scaling the braid group action preserves $\Psi$.   When $\sum
\mu_i = 2$, the signature of $\Psi$ is $(1,n-3)$ by Proposition
\ref{Prop:signature}. It follows that $P\Gamma \subset PU(1,n-3)$.
And so the braid group acts through $P\Gamma$ as automorphisms of
the complex ball $\mathbb{B}^{n-3}$.
\end{proof}

\begin{Cor}\label{Cor:HG-valuedin-Ball}
Assume $\sum_i \mu_i = 2$ and let $|S| = n$.  Then the
multi-valued map $HG_{\mu}: \mathcal{M}_n \rightarrow
\mathbb{B}^{n-3} \subset \mathbb{P}^{n-3}$.
\end{Cor}
\begin{proof}
Because $\sum_i \mu_i = 2$, by Proposition \ref{Prop:signature}
together with Proposition \ref{Prop:IH-Hodge} one sees that
$dim_{\mathbb{C}}(IH_{1,0}(\pone, L_m)) = 1$. Because
$IH_{1,0}(\pone, L_m)$ is positive definite, it follows that the
point $\mathbb{P}(IH_{1,0}(\pone, L_m))$ is an element of the ball
$\mathbb{B}^{n-3} \subset \mathbb{P}(IH_1(\pone, L_m)) =
\mathbb{P}^{n-3}$.  Recall $HG_{\mu}(m)$ is defined to be the
$P\Gamma$ orbit of this point $\mathbb{P}(IH_{1,0}(\pone, L_m))
\in \mathbb{P}^{n-3}$.  By the Proposition,  $P\Gamma$ acts as
automorphisms of the ball, so $HG_{\mu}(m) \subset
\mathbb{B}^{n-3}$.
\end{proof}

\begin{Def}\label{Def:discrete&lattice}
A {\em discrete subgroup} of a Lie group is an infinite subgroup
for which the subspace topology is the discrete topology. A
(group-theoretic) {\em lattice} is a co-finite volume discrete
subgroup of a Lie group. We may sometimes informally refer to
``discrete group" when we really mean a discrete subgroup.
\end{Def}

 \noindent {\bf Example:} Discrete subgroups like
$SL(2, \mathbb{Z})$ and its congruence mod $2$ subgroup
$\Gamma(2)$ are lattices. One can check $\Gamma(2)$ arises as the
monodromy group for the $4$-point case, where $\mu_i =
\frac{1}{2}, \forall i$.

\subsubsection{Monodromy: Braid Action on $IH_1(\pone,
L)$}\label{Subsec:Monodromy}

The monodromy group $\Gamma$ of $\mathcal{L} \rightarrow
\mathcal{P}_n$ is the representation of the spherical $n$-strand
braid group, $\pi_1(\mathcal{P}_n)$, on $IH_1(\pone, L)$.  As a
general reference for standard results on the braid group that we
use, see \cite[Section 5]{HL} and the references contained
therein.

Let $R_{i,i+1}$ denote the braid group ``transposition" element
that braids $s_{i+1}$ about $s_i$ and is the identity on $s_k, k
\neq i+1$. It can be realized by a compactly supported isotopy of
$\pone$ that moves $s_{i+1}$ along a counter-clockwise circle that
encloses $s_i$ and is the identity in neighborhoods of $s_k, k
\neq i+1$. A well-known result is:

\begin{Lemma}\label{Lem:BraidGen}
The spherical braid group on $n$ strands is generated by the
``transpositions" $R_{i, i+1}$ and $R_{n,1}$.
\end{Lemma}

Once a basis for $IH_1(\pone, L)$ is chosen, then the action of
these generators can be written in terms of explicit matrices. For
simplicity we {\em assume} that all of the local monodromies are
non-trivial, that is, $\mu_i \not\in \mathbb{Z}, \, \forall i$. In
Corollary \ref{Cor:IHbasis} we constructed a basis for
$IH_1(\pone, L)$, taken from $\{I_{i,i+1} \equiv \frac{1}{\alpha_i
- 1}\hat{C}_i + \hat{\gamma} +
\frac{1}{1-\alpha_{i+1}}\hat{C}_{i+1} \}_{i \in \{1, \ldots, n-1
\}} \cup \{I_{n,1} \}$. Roughly speaking, any $n-2$ cycles from
this set form a basis. There are two possibilities: either (a) a
point $s_j$ is ``isolated" or (b) some cycle $I_{i,i+1}$ (or
$I_{n,1}$) is ``isolated".

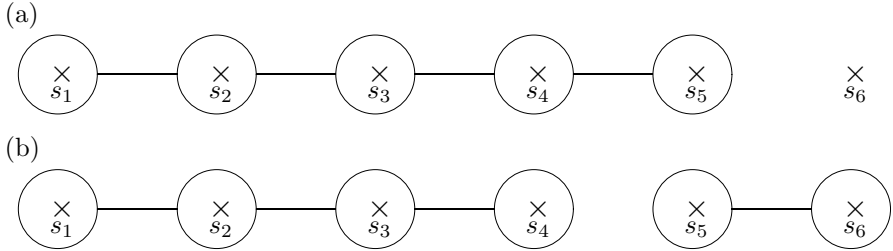
\vspace{.1in}

\begin{picture}(380,50)

\put(0,45){(a)}

\multiput(17,22)(60,0){6}{\large $\times$}

\put(17,16){$s_1$}  \put(77,16){$s_2$} \put(137,16){$s_3$}
\put(197,16){$s_4$} \put(257,16){$s_5$} \put(317,16){$s_6$}


\multiput(20,25)(60,0){5}{\circle{30}}


\multiput(35,25)(60,0){4}{\line(1,0){30}}

\end{picture}

\begin{picture}(380,50)

\put(0,45){(b)}

\multiput(17,22)(60,0){6}{\large $\times$}

\put(17,16){$s_1$}  \put(77,16){$s_2$} \put(137,16){$s_3$}
\put(197,16){$s_4$} \put(257,16){$s_5$} \put(317,16){$s_6$}


\multiput(20,25)(60,0){6}{\circle{30}}


\multiput(35,25)(60,0){3}{\line(1,0){30}}

\put(275,25){\line(1,0){30}}

\end{picture}

\vspace{.1in}

Because all of the points are assumed to have non-trivial
monodromy, the only way to violate the condition of Corollary
\ref{Cor:IHbasis} is with an ``isolated" cycle $I_{i,i+1}$ for
which $\mu_i + \mu_{i+1} \in \mathbb{Z}$.  This will always be a
counter-example to Corollary \ref{Cor:IHbasis}.

\vspace{.1in}

\noindent {\bf Counter-Example:}  Choose a partition $(S_1,S_2)$
which does not satisfy the assumption of Corollary
\ref{Cor:IHbasis}, so that $\sum_{i \in S_1} \mu_i \in \mathbb{Z}$
Then there exists a local system $L_{S_1}$ on $\pone \setminus
S_1$ defined by assigning $\mu_i$ to $s_i \in S$. Then the cycles
$I_{i,i+1}$ cannot be linearly independent, because $IH_1(\pone,
L_{S_1})$ is $(|S_1|-2)$-dimensional and there $(|S_1| - 1)$
cycles.  A simple example is $|S| = 4$, $\mu_i = \frac{1}{2},
\forall s_i \in S$, with $S_1 = \{s_1, s_2 \}$ and $S_2 = \{s_3,
s_4 \}$.

We will give partitions that do not exhibit this pathology a
suggestive name.

\begin{Def}
A partition of $S$ into subsets $S_1$ and $S_2$ where $\sum_{i \in
S_1} \mu_i \not\in \mathbb{Z}$ (or equivalently with $S_2$) {\em
stable partitions}.
\end{Def}

To study $R_{i,i+1}$ it is convenient to take advantage of the
above flexibility in the choice of basis so as to ``isolate" $s_i$
and $s_{i+1}$, like in the above picture (b) of a ``good basis."

\begin{Lemma}\label{Lem:ComplexReflection}
Assume $S_1 = \{s_i, s_{i+1} \}, S_2 = S \setminus S_1$ defines a
stable partition of $S$.  In the good basis above, $R_{i, i+1}$
acts as the identity on the space spanned by the $n-3$ remaining
basis vectors. Furthermore, it acts as an order $k$ complex
rotation, for $k$ the denominator of the fraction (in lowest
terms) $\mu_i + \mu_{i+1}$, on the remaining basis vector $I_{i,
i+1}$.  More specifically, it acts on $I_{i,i+1}$ as
multiplication by $e^{2 \pi \imath (\mu_i + \mu_{i+1})}$.
\end{Lemma}
\begin{proof}
Because $(S_1,S_2)$ is a stable partition, these cycles form a
basis. It is immediate that $R_{i, i+1}$ acts as the identity on
the $n-3$ intersection homology generators associated to $S_2$,
because the isotopy corresponding to the braid action is the
identity away from a small compact set that contains $s_i$ and
$s_{i+1}$ but no other $s_k$.

The action on $I_{i, i+1}$ is more involved.  A formal argument
can be adapted almost {\em mutatis mutandis} from
\cite[Proposition 9.2, pp.46-47]{DM}. Informally it is easy to see
using a ``relative position" argument.  A counter-clockwise motion
of $s_{i+1}$ relative to a fixed $s_i$ may be thought of as a
counter-clockwise motion of $s_i$ relative to a fixed $s_{i+1}$,
with one full loop corresponding to one full loop. The section
therefore is scaled by the local monodromy of each, namely
$\alpha_1 \cdot \alpha_2 = e^{2 \pi \imath(\mu_i + \mu_{i+1})}$.
\end{proof}

\begin{Def}
A finite order complex linear transformation $T$ with a hyperplane
as its fixed point locus is called a {\em complex reflection}. The
{\em mirror} of the reflection is the fixed hyperplane. If $T$
preserves a hyperbolic Hermitian form $\Psi$, then we say $T$ is a
{\em complex hyperbolic reflection}.
\end{Def}

\begin{Prop}\label{Prop:CplxHyperbolicRefl}
Assume the partition $S_1 = \{s_i,s_{i+1} \}, S_2 = S \setminus
S_1$ is a stable partition of $S$. Then $R_{i,i+1}$ is a complex
hyperbolic reflection of order $k$. The mirror of $R_{i,i+1}$ is
the $\Psi$-orthogonal complement of the basis vector $I_{i, i+1}$.
\end{Prop}
\begin{proof}
By Lemma \ref{Lem:ComplexReflection}, $R_{i,i+1}$ is an order $k$
complex reflection.  By Proposition \ref{Prop:PU(1,n)}, $\Gamma$
preserves the hyperbolic Hermitian form $\Psi$ on $IH_1(\pone,
L)$.  By Lemma \ref{Lem:BraidGen} $R_{i,i+1}$ acts on $IH(\pone,
L)$ as a generator of $\Gamma$, and so it must preserve the
hyperbolic structure.

The intersection pairing of $I_{i,i+1}$ with any of the remaining
$n-3$ basis vectors is trivial because their geometric supports do
not intersect.  These vectors associated to $S_2$ therefore span
the $\Psi$-orthogonal complement of $I_{i,i+1}$. By Lemma
\ref{Lem:ComplexReflection}, $R_{i,i+1}$ acts trivially on the
$S_2$ basis vectors, and non-trivially on $I_{i,i+1}$.  Hence the
$\Psi$-orthogonal complement is the mirror of $R_{i,i+1}$.
\end{proof}

\begin{Rmk}
If $S_1 = \{s_i,s_{i+1} \}$ does {\em not} define a stable
partition, then observe that by Lemma
\ref{Lem:Psi-Signature-prelims}, $I_{i,i+1}$ has $\Psi$-length
zero, i.e., is isotropic.
\end{Rmk}

For explicit computations it is useful to have the action of
$R_{i,i+1}$ for all $i$ with respect to a single fixed basis. This
also makes the ring of integers $R = \mathbb{Z}(\zeta_d)$ over
which $\Gamma$ is defined transparent.  Of course, $R$ is the same
as the base ring of the module-theoretic lattice $(IH_1(\pone, L),
\Psi)$, since $\Gamma$ acts as a monodromy group.

\begin{Prop} \label{Prop:GammaReflMatrix}
In the standard basis $I_{j,j+1}, j \in \{1, \ldots, n-2 \}$, the
reflection $R_{i,i+1}$ is a matrix with entries $R_{i,i+1}(a,b)$:
$$ R_{i,i+1}(a,b) = \left\{ \begin{array}{cl}
                        1,  & a = b \neq i \text{ and } |a-b| > 1  \\
                        \alpha_{i} \cdot \alpha_{i+1}, & a = b = i
                        \\
                        1-\alpha_{i+1}, & a = i \text{ and } b =
                        i-1 \\
                        \alpha_{i+1} (1-\alpha_i), & a = i \text{
                        and } b = i+1 \\
                        0, & \text{ elsewhere. }
                    \end{array}  \right.
$$
\end{Prop}

\subsubsection{Uniformization: $INT$ and $\Sigma INT$}

To date we have considered the moduli space $\mathcal{M}_n$ of $n$
{\em distinct} points on $\pone$.  Choosing $\mu$ is equivalent to
choosing a line bundle on $(\pone)^n$, and in fact uniquely
determines a $SL_2(\mathbb{C})$-linearization of the diagonal
$SL_2(\mathbb{C})$ action.  This means there is a well-defined
compact GIT quotient, $\overline{\mathcal{M}}_{n, \mu}$. Let us
denote the quasi-projective stable locus by $\DM$. The key insight
that drives \cite{DM} is that $HG_{\mu}$ extends uniquely over
$\DM$.

The main result of the paper of Deligne and Mostow \cite{DM} is
that, for a finite list of $\mu, HG_{\mu}$ has a single-valued
inverse $\Phi_{\mu}$, and so the bottom map in the following
diagram is an isomorphism of complex analytic spaces.

\vspace{.1in}

\xymatrix{ & & & & & &   & \ \ \mathbb{B}^{n-3} \ar[d]^{\ / \sim} \ar@{->}@(l,u)[1,-1]_{\Phi_{\mu}} \\
& & & & & & {\mathcal{M}_{n,\mu}} \ar@3{->}[ur]^{HG_{\mu}}
\ar[r]^{\cong} & {\Gamma \backslash \mathbb{B}^{n-3}} }

\vspace{.1in}

In fact they show more.  For such $\mu$, the uniformization
extends, as an isomorphism of varieties, to the GIT
compactification $\overline{\mathcal{M}}_{n,\mu}$ (including the
semi-stable points) on the one hand and the Baily-Borel
compactification $\overline{P\Gamma \backslash
\mathbb{B}^{n-3}}^{BB}$ on the other.  In short, ``GIT $=$
Baily-Borel".

Their original sufficiency criterion for $\mu$ is simple to check.

\noindent {\bf Condition $INT$:} Assume that the numbers $\mu_j$
defined by $\alpha_j = e^{2\pi i \mu_j}, 0< \mu_j < 1$ satisfy
$\sum \mu_i = 2$.  For all $ s \neq t$ in $S$ such that $\mu_s +
\mu_t < 1$, require that $(1-\mu_s-\mu_t)^{-1} \in \mathbb{Z}$.

\begin{Thm}[$INT$ \cite{DM}] \label{Thm:INT}
If Condition $INT$ holds, then $\Gamma$ is a lattice in the
projective unitary group $PU(1, n-3)$. Moreover, $\DM \cong
P\Gamma \backslash \mathbb{B}^{n-3}$, and indeed the isomorphism
extends to their GIT and Baily-Borel compactifications as an
isomorphism of varieties.
\end{Thm}

The list of solutions is quite small, and in fact there is only
one solution for $n=7$ and none for $n >7$.  Furthermore, it would
be nice to have a necessary and sufficient condition for
$\Gamma_{\mu}$ to be discrete.  In \cite{Mos1}, Mostow develops a
generalization of $INT$ that largely fulfills that purpose.

\noindent {\bf Condition $\Sigma INT$:} Assume that the numbers
$\mu_j$ defined by $\alpha_j = e^{2\pi i \mu_j}, 0< \mu_j < 1$
satisfy $\sum \mu_i = 2$.  Let $S_1$ be a subset of $S$ with
$\mu_s = \mu_t \, \forall s,t \in S_1$.  For all $s \neq t \in S$
such that $\mu_s + \mu_t < 1$, require that
$$
1 - \mu_s - \mu_t \in \left\{ \begin{array}{cl}
                                \frac{1}{2} \mathbb{Z} & if \ s,t \in
                                S_1 \\
                                \mathbb{Z} & otherwise
                            \end{array}  \right.
$$

\begin{Thm}[$\Sigma INT$ \cite{Mos1}] \label{Thm:SigmaINT}
If Condition $\Sigma INT$ holds then $\Gamma$ is a lattice in
$PU(1, n-3)$. Let $\Sigma$ denote the symmetric group of order
$|S_1|$. Then $\mathcal{M}^{\Sigma}_n \cong P\Gamma_{\Sigma}
\backslash \mathcal{B}^{n-3}$ for a group extension
$\Gamma_{\Sigma}$ of $\Gamma$ by $\Sigma$, and furthermore this
isomorphism extends to their GIT and Baily-Borel compactifications
as an isomorphism of varieties.
\end{Thm}

\begin{Rmk}
The ``$\Sigma$" in $\Sigma INT$ is meant to suggest the symmetric
group. In essence, the idea behind $\Sigma INT$ is to exploit
repeated values in the list $\{ \mu_j \}$ by constructing a
uniformization for $\DM/\Sigma$. So the arguments in the proof
largely reduce to the same arguments used for condition $INT$.
\end{Rmk}

\begin{Rmk}
Whenever an example satisfies $\Sigma INT$, unless otherwise
noted, we by default work with the quotient
$\mathcal{M}^{\Sigma}_n$.
\end{Rmk}

\section{Eisenstein and Gaussian Ancestral Examples}\label{Sec:Ancestral}

\subsection{Automorphisms of lattices}

Let $(L, \Psi)$ be a lattice over $R$, $M$ a sublattice, and $N$
the $\Psi$-orthogonal complement of $M$ in $L$.  Any automorphism
of $L$ restricts to an automorphism $u_M$ of $M$ and an
automorphism $u_N$ of $N$. Conversely, when do the automorphisms
of $M$ extend to automorphisms of $L$?

To address the question we recall some basic ideas from the theory
of lattices.  We assume throughout that $\Psi$ is non-degenerate,
which is automatically true for $\Psi$ as defined in Section
\ref{sec:HGfunctions}.

Let $L^*$ denote the dual lattice $\text{Hom}_R(L, R)$. The form
$\Psi$ induces a map $a_L: L \rightarrow L^*$, given by $x \mapsto
\Psi( \ , x)$.  Since $\Psi$ is nondegenerate, $a_L$ embeds $L$ as
a sublattice of $L^*$ of finite index. (Drawing $L^*$ as the usual
``square" Cartesian lattice, one sees the sublattice $a_L(L)$ is
the standard ``pictorial" representation of the lattice $L$.) Many
of the differences with the theory of vector spaces, where $V^*
\cong V$, are captured by the discrepancy between $L^*$ and
$a_L(L)$.

Let $C(L) := L^*/a_L(L)$, and observe that it is a finite
$R$-module. Furthermore, $\Psi$ determines the Hermitian form
$\Psi^*$ on $L^*$, but now this form is valued in the field of
fractions of $R$, denoted $F(R)$.  This in turn induces a
Hermitian form $\Psi^*_{C(L)}$ on the finite $R$-module $C(L)$.
Note that $C(L)$ is not a ``lattice" {\em per se}, because the
form is valued in the group-theoretic quotient $F(R)/R$.
Nonetheless, it is clear that any unitary transformation of $L$
induces an automorphism of $C(L)$ that preserves $\Psi^*_{C(L)}$.

Assume now that $(\Lambda, \Psi)$ is a {\em unimodular}
$R$-lattice.  To address the question above, consider a {\em
primitive} sublattice $M$ and its $\Psi$-orthogonal complement
$N$.  It turns out that there is a natural isomorphism $\alpha:
C(M) \rightarrow C(N)$ that changes the sign of the forms
$\Psi^*_{C(M)}$ and $\Psi^*_{C(N)}$ but otherwise preserves them.
One can then see:
\begin{Prop} \label{Prop:AutLatticeExtend} \cite[Appendix, pp.43-44]{HL}
Let $M$ be a primitive sublattice of a unimodular lattice $L$.   A
pair of unitary transformations $u_M$ of $M$ and $u_N$ of $N$,
defining a unitary transformation $(u_M,u_N)$ of $M \perp N$, is
an automorphism of $L$ if and only if the following diagram
commutes: \vspace{.1in}

\xymatrix{ & & & & & & C(M) \ar[r]^{\alpha} \ar[d]_{u_M} & C(N) \ar[d]^{u_N} \\
& & & & & & C(M) \ar[r]_{\alpha} & C(N)  }
\end{Prop}

 There are two rings that principally concern us:
the Gaussian and Eisenstein rings of integers.

\begin{Def}\label{Def:GaussEisenRings}
The ring $\mathcal{G}$ of {\em Gaussian} integers is
$\mathbb{Z}[\imath]$, where $\imath = \sqrt{-1}$. The ring
$\mathcal{E}$ of {\em Eisenstein} integers is
$\mathbb{Z}[\omega]$, where $\omega$ is a primitive third root of
unity.
\end{Def}

Remarkably, for $\mathcal{G}$ and $\mathcal{E}$ we don't need an
explicit description of $\alpha$. The results follow from the
properties that $\alpha$ is an isomorphism and (up to sign)
preserves the form.

\begin{Cor}\label{Cor:AutLatticeExtendGaussEisen}
Let $(\Lambda, \Psi)$ be a lattice over the ring $R = \mathcal{G}$
or $\mathcal{E}$.  Let $z \in \Lambda$ be a primitive vector
(i.e., $z$ generates a primitive sublattice $Rz$ in $M$) not of
unit length. Let $\Lambda_0 \subset \Lambda$ be the sublattice
that is $\Psi$-orthogonal to $z$.  Then the map from the
$Aut(\Lambda)$-stabilizer of $\Lambda_0$ to $Aut(\Lambda_0)$ is an
isomorphism.  In particular, any automorphism of $\Lambda_0$
extends uniquely to an automorphism of $\Lambda$.  If $z$ is of
unit length then the ambiguity in the extension is just the
automorphism group of $Rz$, namely $\mathbb{Z}/4\mathbb{Z}$ and
$\mathbb{Z}/6\mathbb{Z}$ for $\mathcal{G}$ and $\mathcal{E}$
respectively.
\end{Cor}
\begin{proof}
The proofs in the Gaussian and Eisenstein cases are analogous.  In
each case the essential point is that $C(Rz)$ is isomorphic to
$R/(r)$ for some $r \in R$.  The unitary automorphisms (those
preserving $\Psi^*_{C(Rz)}$) of $R/(r)$ are easily seen to be one
of these: trivial (if $r$ is a unit), $\mathbb{Z}/6\mathbb{Z}$ if
$R$ is Eisenstein and $r$ not a unit, or $\mathbb{Z}/4\mathbb{Z}$
if $R$ is Gaussian and $r$ not a unit. Consider a unitary
transformation $u_M$ of $\Lambda_0$.  This induces a unitary
transformation of $C(\Lambda_0)$, which, because it is isomorphic
via $\alpha$ to $C(Rz)$, must be an element of the trivial group,
$\mathbb{Z}/6\mathbb{Z}$, or $\mathbb{Z}/4\mathbb{Z}$ according to
the cases above.  Now, $\alpha$ itself acts (accounting for the
sign change) as a unitary automorphism, so it satisfies the same
trichotomy.  In particular, one can always find an automorphism of
$C(Rz)$ to ``undo" $\alpha$ and so make the diagram from the
Proposition commute. The only potential obstruction is that the
requisite automorphism of $C(Rz)$ may not come from an
automorphism $u_N$ of $Rz$.  But the unitary transformations of
$Rz$ are $\mathbb{Z}/6\mathbb{Z}$ for $R$ Eisenstein or
$\mathbb{Z}/4\mathbb{Z}$ for $R$ Gaussian, so in fact one can
always find such a $u_N$, and in particular, as long as $r$ is not
a unit, that $u_N$ is determined uniquely by $u_M$. In other
words, any automorphism of $\Lambda_0$ extends uniquely to an
automorphism of $\Lambda$.  If $r$ is a unit, then the ambiguity
is precisely the group of units in $\mathcal{G}$ or $\mathcal{E}$.
\end{proof}

In general, given a locally symmetric space (here, a ball
quotient), it can be quite difficult to identify locally symmetric
subspaces (here, subball quotients).

\begin{Def}\label{Def:subballquot}
Let $P\Gamma$ be a discrete subgroup of $PU(1,n)$.  Let
$\mathbb{B}^k$ be a subball of $\mathbb{B}^n$.  In particular,
$\mathbb{B}^k$ is cut out by a (projective) linear constraint on
the ambient $\mathbb{P}^n$.  Let $P\Gamma_{Stab}$ denote the
subgroup of $P\Gamma$ that preserves $\mathbb{B}^k$.  Consider the
image of $\mathbb{B}^k$ in the ball quotient $P\Gamma \backslash
\mathbb{B}^n$.  We say that this image is a subball quotient if
the map factors through an inclusion of $P\Gamma_{Stab} \backslash
\mathbb{B}^k$ in $P\Gamma \backslash \mathbb{B}^n$.

\vspace{.1in}

\xymatrix{ & & & & & {\mathbb{B}^k} \ar[d] \ar[r] & {\mathbb{B}^n} \ar[d] \\
& & & & & {P\Gamma_{Stab} \backslash \mathbb{B}^k} \ar[r] &
{P\Gamma \backslash \mathbb{B}^n} }
\end{Def}

The above Corollary tells us that, for a unimodular lattice over
the Gaussian or Eisenstein integers, ``arithmetically-defined"
hyperballs $\mathbb{B}^{n-1}$ correspond to codimension 1 subball
quotients.  Induction yields:

\begin{Cor}\label{Cor:GaussEisenSubballQuot}
For $\Lambda$ a unimodular lattice of hyperbolic signature over $R
= \mathcal{G}$ or $\mathcal{E}$, any primitive hyperbolic
sublattice $\Lambda_0$ defines a subball quotient:
$$PAut(\Lambda_0) \backslash \mathbb{B}^{n-1} \subset PAut(\Lambda)
\backslash \mathbb{B}^n
$$
\end{Cor}
\begin{proof}
When $\Lambda_0$ is the $\Psi$-orthogonal complement of a
primitive vector in $\Lambda$ and is of hyperbolic signature this
is a restatement of the previous Corollary. It is clear that the
intersection of subball quotients is again a subball quotient. So,
because $\Lambda_0$ is the $\Psi$-orthogonal complement of some
primitive lattice, the statement follows by induction.
\end{proof}

\subsection{Organizing principle: Descendants by collision}

All the $\Gamma$ discussed in this section are assumed to satisfy
$\Sigma INT$ (and so in particular are group-theoretic lattices,
i.e., discrete subgroups of $PU(1,n-3)$, defined over some ring of
integers), thus $\DM \cong P\Gamma \backslash \mathbb{B}^{n-3}$. A
collision between two points $s_i$ and $s_j$ is identified with
the complement of a lattice vector, yielding a codimension 1
subball quotient in $P\Gamma \backslash \mathbb{B}^{n-3}$.  This
is implicit in Deligne and Mostow's main theorems, as it is a part
of the extension of the uniformization over the stable boundary of
\DM. To be explicit, using the notation introduced in Section
\ref{Subsec:Monodromy}:

\begin{Lemma}\label{Lem:CollisionMirror}
Assume $\{s_i,s_j\} \cup S \setminus \{s_i,s_j\}$ is a stable
partition.  Let $S_{i,j}$ denote the sublocus consisting of all
configurations of points for which $s_i$ and  $s_j$ share a
coordinate (i.e., have ``collided").  The image of the principal
branch of $HG_{\mu}$ restricted to $S_{i,j}$ is the mirror of
$R_{i,j}$. Equivalently it is the $\Psi$-orthogonal complement of
the vector in $\mathbb{B}^{n-3}$ assigned to $I_{i,j}$.
\end{Lemma}
\begin{proof}
Because it is a stable collision, $S_{i,j}$ is a nonempty subset
of $\DM$. By Theorems \ref{Thm:INT} and \ref{Thm:SigmaINT},
$HG_{\mu}$ is well-defined on $S_{i,j}$. For convenience, relabel
the points to be $s_i$ and $s_{i+1}$. $HG_{\mu}(m)$ is valued in
$\mathbb{P}(IH_1(\pone, L_m))$.  Consider the good basis that
``isolates" $I_{i,i+1}$, which exists by Corollary
\ref{Cor:IHbasis} because this is a stable partition of $S$.
$HG_{\mu}(m)$ is (the projective image of) a linear combination of
these basis vectors, or equivalently by Proposition
\ref{Prop:CplxHyperbolicRefl}, of $I_{i,i+1}$ and the basis
vectors in its $\Psi$-orthogonal complement. When $s_i(m) =
s_{i+1}(m)$ via a path $0$ to $m$ that does not cross a branch, a
good basis for $IH_1(\pone, L_m)$ is precisely (the flat translate
of) the basis for the $\Psi$-orthogonal complement of $I_{i,
i+1}$, denoted $I_{i,i+1}^{\perp}$; that is, $IH_1(\pone, L_m)$ is
the mirror of the complex reflection $R_{i,i+1}$.

In particular, for such $m, IH_{1,0}(\pone,L) \subset
I_{i,i+1}^{\perp}$, or equivalently, $HG_{\mu}(m) \in
\mathbb{P}(I_{i,i+1}^{\perp}) \cap \mathbb{B}^{n-3}$. Because
$HG_{\mu}$ is by assumption a uniformization, the image of
$S_{i,i+1}$ is an open subset of the hyperball
$\mathbb{I_{i,i+1}^{\perp}} \subset \mathbb{B}^{n-3}$. But the
sublocus is also a closed subset, so again because $HG_{\mu}$ is
an isomorphism with $P\Gamma \backslash \mathbb{B}^{n-3}$, it must
map to a closed set and hence the image is both open and closed in
the hyperball, and so is the full hyperball.
\end{proof}

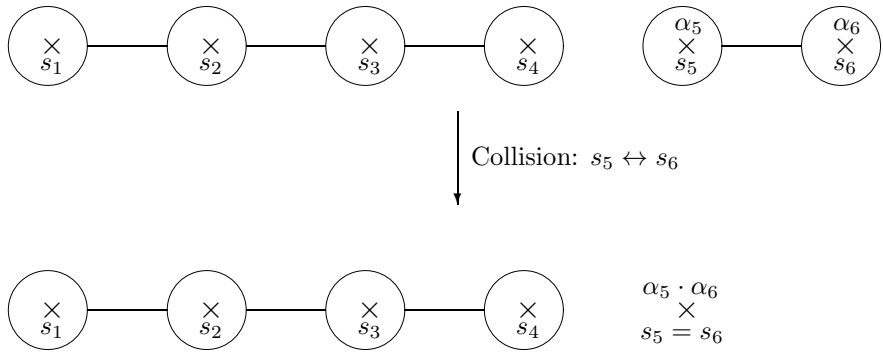
\vspace{.1in}
\begin{picture}(400,150)

\multiput(17,122)(60,0){6}{\large $\times$}

\put(17,116){$s_1$}  \put(77,116){$s_2$} \put(137,116){$s_3$}
\put(197,116){$s_4$} \put(257,116){$s_5$} \put(317,116){$s_6$}

\put(257,131){$\alpha_5$}  \put(317,131){$\alpha_6$}


\multiput(20,125)(60,0){6}{\circle{30}}


\multiput(35,125)(60,0){3}{\line(1,0){30}}

\put(275,125){\line(1,0){30}}


\put(175,100){\vector(0,-1){35}}

\put(180,80){Collision: $s_5 \leftrightarrow s_6$}


\multiput(17,22)(60,0){5}{\large $\times$}

\put(17,16){$s_1$}  \put(77,16){$s_2$} \put(137,16){$s_3$}
\put(197,16){$s_4$}  \put(244,14){$s_5 = s_6$}

\put(245,31){$\alpha_5 \cdot \alpha_6$}


\multiput(20,25)(60,0){4}{\circle{30}}


\multiput(35,25)(60,0){3}{\line(1,0){30}}

\end{picture}

\vspace{.1in}

\begin{Rmk}
By induction, a general collision sublocus is just the
intersection of mirrors, and so the orthogonal complement of a
collection of vectors.
\end{Rmk}

There are precisely four Deligne-Mostow lattices that are
generated by $\mu$ with all the $\mu_i$ equal-valued.  These are
the equally weighted $n=4, 5,6,8,$ and $12$ point examples.  For
$n=6$ and $12$, $R$ is the Eisenstein ring $\mathbb{Z}[\omega]$
(here $\omega$ is a primitive sixth root of unity), whereas for
$n=4$ and $8$, $R$ is the Gaussian ring $\mathbb{Z}[\imath]$.

One special feature of these lattices is that $\Gamma_{\Sigma}$ is
the {\em full} automorphism group of the associated
module-theoretic lattice, rather than just a subgroup. In these
two cases results can be found, from a different perspective and
in different language, in the recent literature.

\begin{Thm}[\cite{MY}] \label{Thm:GaussianAncestral}
 The uniformizing group $\Gamma_G$ for the Gaussian ancestral example, that is, the
Deligne-Mostow example for $8$ equally weighted $(\mu_i =
\frac{1}{4})$ points, equals the full automorphism group of the
corresponding lattice. That is, $\Gamma_G =
Aut(\mathbb{Z}[\imath], \Psi_G)$.
\end{Thm}

\begin{Thm}[\cite{All2}] \label{Thm:EisensteinAncestral}
 The uniformizing group $\Gamma_E$ for the Eisenstein
ancestral example, that is, for $12$ equally weighted $(\mu_i =
\frac{1}{6})$ points, equals the full automorphism group of the
corresponding lattice. That is, $\Gamma_E =
Aut(\mathbb{Z}[\omega], \Psi_E)$.
\end{Thm}

\begin{Def}\label{Def:Ancestral}
We call these the {\em ancestral} Deligne-Mostow lattices. The $8$
point case we call the {\em Gaussian ancestral lattice} and
similarly the $12$ point case we call the {\em Eisenstein
ancestral lattice}.  By a {\em descendant} lattice, we mean the
subgroup of an ancestral lattice which is the stabilizer subgroup
for a subconfiguration space (collision sublocus).
\end{Def}

\begin{Thm}\label{Thm:DescendantAutofLattices}
Descendants of the  ancestral lattices are themselves automorphism
groups of the (module-theoretic) sub-lattices.  The collision loci
in $DM(1^8)$ and $DM(1^{12})$ are orbifold subball quotients.
\end{Thm}
\begin{proof}
It suffices to check for codimension $1$, the rest follow by
induction.  By Lemma \ref{Lem:CollisionMirror} the image under (a
branch of) $HG_{\mu}$ of a stable collision of a pair of points is
the $\Psi$-orthogonal complement of a vector $I_{i,j}$.  The
vector lies on the lattice $\Lambda$ in $IH_1(\pone, L)$, so the
complement defines a sublattice $\Lambda_0$. By Lemma
\ref{Lem:Psi-Signature-prelims}, $I_{i,j}$ has negative length, so
its $\Psi$-orthogonal complement is hyperbolic.   By Corollary
\ref{Cor:AutLatticeExtendGaussEisen}, any automorphism of
$\Lambda_0$ therefore extends to an automorphism of $\Lambda$.  So
the stabilizer subgroups are in fact themselves automorphism
groups of sub-lattices.  The non-uniqueness of the extension is
the order of $R_{i,j}$ as a complex reflection, which is
non-trivial by Lemma \ref{Lem:ComplexReflection}, so these are
orbifold loci.
\end{proof}

\begin{Rmk}
Making use of the three common meanings of ``lattice" in
mathematics --- poset, group theoretic, and module theoretic ---
this Theorem tells us we have described, amusingly, a ``lattice of
lattices which are automorphisms of lattices".
\end{Rmk}

It is straightforward to observe that the equally weighted $n=6$
case, defined by $\mu = (\frac{1}{3}, \ldots, \frac{1}{3})$, is a
descendant of the Eisenstein example, where the $12$ points have
all collided in pairs.  Since the $n=5$ case is two complex
dimensional, it has no descendants of Deligne-Mostow type (and
precisely one descendant of dimension $1$).

\begin{Cor} \label{Cor:TwoAncestral}
The only equally weighted examples with proper Deligne-Mostow
descendants are the Gaussian and Eisenstein ancestral examples.
\end{Cor}

\begin{Cor}\label{Cor:DescendantListObserve}
For $n > 7$, all the Deligne-Mostow lattices are (finite index
sublattices of) descendants of the Eisenstein and Gaussian
ancestral lattices. Similarly, all but one of the $n = 7$ examples
is a descendant, and a number of the remaining ones ($n = 5, 6$)
are as well (see Chart). Furthermore, the original list due to
Mostow is in error, for it misses a number of descendant
solutions.
\end{Cor}
\begin{proof}
This follows by direct observation and comparison with Mostow's
chart \cite{Mos1}.
\end{proof}
\begin{Rmk} \label{Rmk:Thurston}
Thurston, working on the problem of enumerating flat metrics with
cone singularities on $S^2$, corrected Mostow's computations by a
computer check, and his list should be complete \cite{Thur}.
\end{Rmk}

\begin{Rmk}
We show in \cite{Dor1} that the moduli space of cubic surfaces
inherits a ball quotient structure, agreeing with that discovered
in \cite{DKvG}, from the Eisenstein descendant $DM(2^5,1^2)$,
which is one of the examples missed by Mostow's tables.
\end{Rmk}

\section{Pull-back Construction}\label{sec:Pullback}

\subsection{Intersection homology under pull-back}\label{SubSec:PullbackIH}

Fix a finite subset $T \subset \pone$ and define a rank $1$
Deligne-Mostow local system $l_T$ on $\pone \setminus T$ with
monodromy $\nu$ and ring of definition $R$.  Consider a map $\pi:
\pone \rightarrow \pone$. Denote the inverse image sheaf, known
henceforth as the {\em pull-back local system},  on $\pone
\setminus \pi^{-1}(T)$ by $\pi^*l_T$.  Because it is rank $1$,
$\pi^*l_T$ is determined by local monodromies at the elements of
$\pi^{-1}(T)$ (by Proposition \ref{Prop:rk1localmonod}), which in
turn can be expressed in terms of $\nu$ and the ramification
indices of $\pi$. More precisely:

\begin{Lemma}\label{Lem:PullbackDM-Lsyst}
 Let $p_{i,j}$ denote the points of the set $\pi^{-1}(t_j)$ and let $r_{i,j}$ denote the
ramification index of $\pi$ at $p_{i,j}$.  Then $\pi^*l_T$ is the
Deligne-Mostow local system on $\pone \setminus \pi^{-1}(T)$
defined by the local monodromy data $r_{i,j} \cdot \nu_j$ at
$p_{i,j}$.  It contains the pull-back local subsystem
$\pi^*l_T(R)$ with fiber $R$.
\end{Lemma}

We now study how $\pi$ induces maps on intersection homology.  One
approach is, using the formalism due to Deligne developed in
\cite{GM2}, to define $\pi_*$ and its adjoint map $\pi^*$ at the
level of the intersection chain complexes for the cover $X$ and
the base $Y$.  To avoid introducing new notation, it is more
direct to follow \cite{GMLef}, and use the following definition.

\begin{Def}\label{Def:PlacidMap}
A subanalytic map $f: X \rightarrow Y$ between two subanalytic
pseudo-manifolds is called {\em placid} if there exists a
subanalytic stratification of $Y$ such that for each stratum $S$
in $Y$ we have
$$ codim_X f^{-1}(S) \geq codim_Y (S) $$
\end{Def}

Any branched covering is placid, so in particular $\pi$ is placid,
where the strata for $Y$ are given by $(T, \pone \setminus T)$ and
those for $X$ by $(\pi^{-1}(T), \pone \setminus \pi^{-1}(T))$.

Intersection homology is a bivariant functor for placid maps,
where the contravariant induced map may shift degrees.  Although
the following Proposition is proven in \cite[Proposition
4.1]{GMLef} for intersection homology valued in the trivial rank
$1$ rational local system (i.e., the constant sheaf with stalk
$\mathbb{Q}$), its proof immediately generalizes to intersection
homology valued in a rank $1$ local system $L \rightarrow Y$ and
in the pull-back local system $f^*L \rightarrow X$.
(Alternatively, one can prove this formally, for the topological
definition of placid maps, using Deligne's construction of
intersection homology.)

\begin{Prop}\label{Prop:Placid-Pi-IH}
Suppose $f: X \rightarrow Y$ is a placid map.  Let $L \rightarrow
Y$ be a rank $1$ local system, and let $f^*L$ denote the pull-back
local system on $X$.  Then pushforward of chains and pull-back of
generic chains induces homomorphisms on intersection homology,
\begin{eqnarray*}
f_*: & IH_i(X, f^*L) & \rightarrow IH_i(Y,L) \\
f^*: & IH_i(Y, L) & \rightarrow IH_{i + dim(X) - dim(Y)}(X, f^*L)
\ .
\end{eqnarray*}
In particular, $\pi: \pone \rightarrow \pone$ induces a map
$\pi^*: IH_1(\pone, l_T) \rightarrow IH_1(\pone, \pi^*l_T)$.  This
map respects the intersection pairing.
\end{Prop}
\begin{Rmk}\label{Rmk:Placid-Pi-IH-Ring} Indeed, one
does not need the fiber to be a field; as is remarked in
\cite{GMLef} after the proof of the Proposition, the same argument
carries over for any coefficient ring $R$.  The same result thus
holds for any local subsystems $L(R)$, with $R$ a subring of
$\mathbb{C}$. In that event, $\pi^*$ is a map of $R$-modules.
\end{Rmk}

Furthermore, (using the differential form model for intersection
cohomology) pulling back a (anti-)holomorphic form via an
algebraic map yields a (anti-)holomorphic form, so the orthogonal
decomposition into $IH^{1,0} \oplus IH^{0,1}$ is respected.  The
isomorphism with intersection homology via the intersection
pairing,  $IH^1(\pone, L) \rightarrow IH_1(\pone, L), z \mapsto (
\ , z)$, tautologically respects the orthogonal decomposition.
Thus the map $\pi^*$ on intersection homology also respects the
Hodge decomposition.

\begin{Prop} \label{Prop:Pi-IH-Hodge}
The map $\pi^*: IH_1(\pone, l_T) \rightarrow IH_1(\pone,
\pi^*l_T)$ preserves the intersection pairing and hence the
Hermitian form $\Psi$, in the sense $\Psi(\alpha, \beta) =
\Psi(\pi^*(\alpha), \pi^*(\beta))$.  In addition $\pi^*$ respects
the orthogonal direct sum (Hodge) decomposition, so that the
subspace $\pi^*(IH_{1,0}(\pone, l_T)) \subset IH_{1,0}(\pone,
\pi^*l_T)$ and $\pi^*(IH_{0,1}(\pone, l_T)) \subset
IH_{0,1}(\pone, \pi^*l_T)$.
\end{Prop}
\begin{Rmk}
This follows more formally from work of Saito on mixed Hodge
modules.  See \cite{Saito1}. Furthermore, as a consequence, one
can amplify Remark \ref{Rmk:Placid-Pi-IH-Ring}.  If one works with
local subsystems whose fibers are the ring $R \subset
\mathbb{C}^*$, then $\pi^*$ is a map of Hermitian lattices over
$R$.
\end{Rmk}

\subsection{Hurwitz spaces and \Spi}

Now we vary $\pi$ while preserving the ramification behavior over
the fixed branch locus $T$.  For any $\pi$ in this family, the
pull-back local system $\pi^*l_T$ will have the same monodromy
data $\mu$.  As $\pi$ varies, the coordinates of the points of
$\pi^{-1}(T)$ vary.

\begin{Def} \label{Def:S-subset-PiInverse(T)}
Let $S \subset \pi^{-1}(T)$ denote the subset of points with
nontrivial local monodromy in $\pi^*l_T$.
\end{Def}
In particular $S$ varies with $\pi$; we write this dependence as
$S(\pi)$. Let us be more precise:
\begin{Def}\label{Def:T-ramification}
By the {\em $T$-ramification class} of $\pi$, denoted
$\mathcal{H}_{\pi}$, we mean the subset of all maps $\pi^{\prime}:
\pone \rightarrow \pone$ satisfying three conditions.

\begin{enumerate}

\item $\pi^{\prime}$ has the same degree as $\pi$.

\item The ramification indices $r_{i,j}$ over points $t_j \in T$
are the same for $\pi^{\prime}$ and $\pi$.

\item $\pi^{\prime}$ is in the same connected component as $\pi$
(with respect to the subspace topology of the standard topology on
the space of maps between compact sets).
\end{enumerate}
\end{Def}
The notation is meant to emphasize the link with Hurwitz spaces,
i.e., spaces of curve covers up to equivalence, since these
self-maps of $\pone$ are curve covers with constrained
ramification.  Equivalence of curve covers is given by the
$PGL_2(\mathbb{C})$ action on $\pone = \mathbb{P}(V)$, which lifts
to an $SL_2(\mathbb{C})$ action on $V$.  Of course,
$SL_2(\mathbb{C})$ also acts on the sets $\pi^{-1}(T)$ and $S$ via
$Sym^{|\pi^{-1}(T)|}(V)$ and $Sym^{|S|}(V)$ respectively. We
denote the induced $SL_2(\mathbb{C})$-equivariant algebraic maps
by
\begin{eqnarray*}
\mathcal{H}_{\pi}  \stackrel{i_{T}}{\longrightarrow}
\mathbb{P}^{|\pi^{-1}(T)|}  \stackrel{p_S}{\longrightarrow}  \mathbb{P}^{|S|} \\
\pi^{\prime}   \mapsto  {\pi^{\prime}}^{-1}(T)  \mapsto
S(\pi^{\prime}).
\end{eqnarray*}

We are interested in the space of all configurations
$S(\pi^{\prime}) \subset \pone$ where $\pi^{\prime} \in
\mathcal{H}_{\pi}$.
\begin{Def}
Let $\mathcal{S}_{\pi}$ denote the $SL_2(\mathbb{C})$-quotient of
the image subvariety $p_s \circ \imath_T(\mathcal{H}_{\pi})
\subset \mathbb{P}^{|S|}$.
\end{Def}

A curve cover $\pi \in \mathcal{H}_{\pi}$ is simply a rational
function on $\pone$. Let $V \cong \mathbb{C}^2$ with coordinates
$(u,v)$. The set of all degree $d$ maps $\pi: \pone \rightarrow
\pone$ is given by:
$$
\big\{ \frac{N(u,v)}{D(u,v)} \big\vert N,D \in Sym^d(V) \big\} =
\mathbb{P}(Sym^d(V) \oplus Sym^d(V)) \ .
$$
In particular $\mathcal{H}_{\pi}$ is an
$SL_2(\mathbb{C})$-invariant subvariety of $\mathbb{P}(Sym^d(V)
\oplus Sym^d(V))$.

\begin{Prop}
$\imath_T$ is injective
\end{Prop}
\begin{proof}
Observe that the numerator $N(u,v)$ determines $\pi^{-1}(0)$ and
the denominator $D(u,v)$ determines $\pi^{-1}(\infty)$. Conversely
these two sets of points determine $\pi$ up to scaling.

Since $|T| \geq 3$, use an automorphism of $\mathbb{P}^1$ to
assign $t_0 = 0$, $t_1 = 1$, and $t_2 = \infty$. Then
$\pi^{-1}(t_0)$ and $\pi^{-1}(t_2)$ determine $\pi$ up to scaling.
But if $\pi^{\prime} = \lambda \pi$, then
$\pi^{\prime}(\pi^{-1}(t_1)) = \lambda$, so in fact
$\pi^{-1}(t_1)$ determines the scaling factor.
\end{proof}

The dimension of $\mathcal{H}_{\pi}$ is easy to compute, as it is
essentially an application of Riemann-Hurwitz.
\begin{Prop}\label{Prop:T-ram-dimcount}
When it exists, the $T$-ramification class $\mathcal{H}_{\pi}$,
where $\pi$ is degree $d$, is a $SL_2(\mathbb{C})$-invariant
subvariety of $\mathbb{P}(Sym^d(V) \oplus Sym^d(V))$, with
codimension equal to $\sum_{i,j} (r_{i,j} - 1)$.
\end{Prop}
\begin{proof}
The Riemann-Hurwitz formula here states $2(d-1) = \sum_k (r_k-1)$
where the sum is over {\em all} ramification points.  Up to
$SL_2(\mathbb{C})$-equivalence, generically the set of covers is
in one-to-one correspondence with the set of coordinates of the
ramification points. The requirement that a ramification point
with index $r_{i,j}$ map to a specific $t_j \in \pone$ is
therefore a codimension $r_{i,j}-1$ condition, and all of these
are independent.
\end{proof}
\begin{Rmk}
Equivalently, the dimension, accounting for
$SL_2(\mathbb{C})$-equivalence, is the number of ``free'' or
``excess'' simple (order 2) ramification points (i.e., those not
in $\pi^{-1}(T)$) allowed by Riemann-Hurwitz.
\end{Rmk}

\subsection{Restricting hypergeometric functions: $\Spi$ and subball
quotients}

Let $n = |S|$, and let the monodromy data for $S$ be $\mu$. If
$(S, \mu)$ can be realized via pullback by $\pi$, then $\Spi
\subset \DM$.  We consider the multi-valued hypergeometric
function $HG_{\mu}$, defined on $\DM$, restricted to $\Spi$. The
restricted hypergeometric function satisfies a linear constraint.
More precisely:

\begin{Lemma}\label{Cor:Spi-HGLinearConstraint}
$\omega_{\mu}(\Spi)$ is $\Psi$-orthogonal to the well-defined
marked subspace $\pi^*(IH_{0,1}(\pone, l_T))$. In particular, a
branch of $HG_{\mu}(\Spi)$ lies in a subball $\mathbb{B}^k \subset
\mathbb{P}((\pi^*(IH_{0,1}(\pone, l_T)))^{\perp})$.
\end{Lemma}
\begin{proof}
Let $R$ be the ring of integers in $\mathbb{Q}(\zeta_d)$. The
restriction of $\mathcal{L}_{\mu} \rightarrow \mathcal{M}_{|S|}$
to $\Spi$ is a local system with a marked local subsystem
$\pi^*(IH_{0,1}(\pone, l_T(R)))$.  The fibers are of the subsystem
are:
$$\pi^*(IH_{0,1}(\pone, l_T(R))) \subset IH_1(\pone,\pi^*l_T(R)) \cap
IH_{0,1}(\pone, \pi^*l_T).$$
 The marked subspace
$\pi^*(IH_{0,1}(\pone, l_T))$ is independent of choice of $\pi$ in
$\Spi$, because $\Spi$ is connected and because $\pi^*(IH_1(\pone,
l_T(R)))$ is a sublattice and so is invariant under continuous
deformations of $\pi$. Then, for any $\pi \in \Spi$,
$\omega_{\mu}(S(\pi)) \in (\pi^*(IH_{0,1}(\pone, l_T)))^{\perp}$
by Proposition \ref{Prop:Pi-IH-Hodge}. Choosing consistent
coordinates by not extending over branch loci, the rest follows
immediately from the definition of $HG_{\mu}$.
\end{proof}

So a branch of the restricted hypergeometric function always lies
in a subball.  We use this fact, applied to Deligne-Mostow
uniformizations of Eisenstein or Gaussian type, to produce $\Spi$
that give subball quotients.  Let $\overline{\Spi}$ denote the
closure of $\Spi$ in $\DM$.

\begin{Thm}\label{Thm:MyMain}
Let $\mu$ be of Eisenstein or Gaussian type satisfying $\Sigma
INT$. If $dim_{\mathbb{C}}(\Spi) =
dim_{\mathbb{C}}(\pi^*(IH_{0,1}(\pone, l_T(R))))$, then
$\overline{\Spi} \subset \DM$ is the Baily-Borel compactification
of the subball quotient $\Gamma_{Stab} \backslash
\mathbb{B}((\pi^*(IH_{0,1}(\pone, l_T)))^{\perp})$.
\end{Thm}
\begin{proof}
We know $\overline{\Spi}$ is an algebraic subvariety of
$\mathcal{M}_{n,\mu}$.  Furthermore, the subball quotient
$\overline{P\Gamma_{Stab} \backslash \mathbb{B}^k}^{BB}$ is an
irreducible algebraic subvariety of $\overline{P\Gamma \backslash
\mathbb{B}^{n-3}}$.  Because $\Phi$ is an isomorphism, and
$\Phi(\Spi) \subset P\Gamma_{Stab} \backslash
\mathbb{B}((\pi^*(IH_1(\pone, l_T)))^{\perp})$, we see $\Spi$ and
its closure inject as subvarieties.  Since the only
equal-dimensional closed subvariety of an irreducible variety is
the irreducible variety itself, as long as the dimensions are
equal one concludes $\Phi$ restricts to give an isomorphism of
these two varieties.
\end{proof}

It is therefore important to compute the dimension of $\Spi$. We
know it equals the dimension of $\mathcal{H}_{\pi}$ precisely when
the map $p_S$ restricted to $\mathcal{H}_{\pi}$ is generically
finite-to-one. One interesting class of examples is when all the
non-trivial monodromy lies over a single point of $T$.

\begin{Lemma}\label{Lem:dimHpiSpi}
Assume $S = \pi^{-1}(t_j)$ for some $t_j \in T$.  If, for some
$k$, $t_k$ has local monodromy $\nu_k = \frac{n_k}{d_k}$ such that
$d_k > 2$, then $dim_{\mathbb{C}}(\Spi) =
dim_{\mathbb{C}}(\mathcal{H}_{\pi})$.
\end{Lemma}
\begin{proof}
Using the $SL_2(\mathbb{C})$-action one may equivalently assume $S
= \pi^{-1}(\infty)$, that $t_k = 0$, and that some $t_i = 1$. We
claim there are only finitely many $\pi^{\prime}$ in
$p_S^{-1}(S)$. For ${\pi^{\prime}}^{-1}(0)$ to consist of points
with trivial local monodromy, each point must be ramified to order
a multiple of at least 3 which is a codimension at least
$\frac{d}{3}(3-1)$ condition. Similarly and independently, over
$1$ every point is ramified to order a multiple of at least 2,
which is codimension at least $\frac{d}{2}(2-1)$. Each fiber of
$p_S$ must therefore be less than $d$ dimensional.  But the
denominator of $\pi$ is determined up to scaling by $S$, which is
a $d$-dimensional condition.  So the generic fiber of $p_S$ is
finite-to-one.
\end{proof}

\subsection{Key examples: $|T| = 3$ and moduli spaces of inhomogeneous
binary forms}

We now explicitly work out the simplest examples.  Assume $|T|=3$,
$\Sigma \mu_i = 2$, $\Sigma \nu_i = 1$, $\mu$ satisfies $\Sigma
INT$, and $\mu$ is Eisenstein or Gaussian.  Specializing our
previous results, we obtain:

\begin{Cor}\label{Cor:Codim1IsBallQuot}
$\overline{\Spi} \subset \DM$ is a subball quotient if and only if
$\pi^*$ is non-trivial and $\Spi$ is codimension $1$.
\end{Cor}
\begin{proof}
Here $IH_1(\pone \setminus T, l_T)$ is one-dimensional and purely
anti-holomorphic.  The image under pull-back is either trivial or
one-dimensional, and purely anti-holomorphic.
\end{proof}

It is worthwhile to completely classify the solutions in a special
case.  By Corollary \ref{Cor:Codim1IsBallQuot}, we need to compute
$dim_{\mathbb{C}}(\Spi)$.  Lemma \ref{Lem:dimHpiSpi} suggest we
consider $S = \pi^{-1}(t_i)$. Restrict further to the case where
all points over a given $t_j \in T$ have the same ramification
index. We think of this as a weak form of a ``Galois" condition on
$\pi$, and so define:
\begin{Def}\label{Def:PropertyG}
The pair $(\pi, T)$ possesses {\em property $G$} if, for $t_j \in
T$, all the points in $\pi^{-1}(t_j)$ have the same ramification
index $r_{t_j}$.
\end{Def}
Up to automorphisms of $\pone$, we may take $T = \{ 0, 1, \infty
\}$.
\begin{Prop}\label{Prop:Pi-PropertyG}
Let $T = \{ 0, 1, \infty \}$ and assume $\pi$ is degree $d$. Then
$\pi$ has property $G$, with $r_0 = a$, and $r_1 = b$ if and only
if
$$
\pi = \frac{A^a(u,v)}{A^a(u,v)+B^b(u,v)}, \ deg(A) = \frac{d}{a},
\  deg(B)= \frac{d}{b} \ ,
$$  where neither $A(u,v)$ nor $B(u,v)$ have repeated roots,
and where $A^a(u,v) + B^b(u,v)$ either has no repeated roots or is
of the form $C(u,v)^c$ for $c | d$ where $C(u,v)$ has no repeated
roots.
\end{Prop}
\begin{proof}
The points in $\pi^{-1}(x)$ are the solutions to $\pi(u,v) = x$.
The ramification index of a point is its multiplicity as a
solution.   Write $\pi = N(u,v)/D(u,v)$.  The numerator $N(u,v)$
and the denominator $D(u,v)$ are both degree $d$ homogeneous
polynomials.

Let $\pi$ satisfy the assumptions.  The points of $\pi^{-1}(0)$
are simply the roots of $N(u,v)$.  Property $G$ says they all must
have the same multiplicity, and so there are $k_1 = \frac{d}{a}$
distinct roots each with multiplicity $a$. Therefore $N(u,v) =
A^a(u,v)$. Similarly, $\pi^{-1}(1)$ consists of the roots of
$N(u,v)-D(u,v)$. By the assumptions on $\pi$ this must have $k_2 =
\frac{d}{b}$ distinct roots each of multiplicity $b$ and so equals
$B^b(u,v)$, implying $D(u,v) = A^a(u,v) + B^b(u,v)$. Finally, the
roots of $D(u,v)$ are the points of $\pi^{-1}(\infty)$ and so by
property $G$ must have equal multiplicities; hence $D(u,v) =
C^c(u,v)$ for some $c | d$.

Conversely, given the explicit form for $\pi$, successively set
$\pi$ equal to $0, 1,$ and $\infty$, and solve.  By assumption
$A(u,v)$ and $B(u,v)$ have no repeated roots.  Therefore
$\pi^{-1}(0)$ is a set of $k_1 = \frac{d}{a}$ points of
ramification index $a$, and $\pi^{-1}(1)$ is a set of $k_2 =
\frac{d}{b}$ points of ramification index $b$.  Likewise, the
assumption that $D(u,v) = A^a(u,v) + B^b(u,v) = C^c(u,v)$
guarantees that the points of $\pi^{-1}(\infty)$ all have
ramification index $c$ (possibly equal to $1$).  Therefore by
definition $\pi$ has property $G$.
\end{proof}

We want to enumerate all $\nu$ such that $(\nu, \mu, \pi)$ satisfy
all of our operating assumptions.

\vspace{.1in}

\noindent  {\bf Summary of Assumptions:}
\begin{enumerate}
\item $T = \{ 0, 1, \infty \}$ (arranged by automorphisms of
$\pone$)

\item $\nu$, the monodromy data defining the DM local system $l_T
\rightarrow \pone \setminus T$, satisfies $\sum \nu_j = 1$

\item $\pi$ satisfies property $G$

\item $\mu$ is the monodromy data for the pull-back local system
$\pi^*l_T$, such that:

    \begin{enumerate}

        \item $S = \pi^{-1}(\infty)$, that is, $\mu_i \not\in
        \mathbb{Z}$ precisely for $s_i \in \pi^{-1}(\infty)$

        \item the sum of the non-integral $\mu_i$ equals $2$

        \item Let $m$ be the lowest common denominator of the
        $\mu_i$.  The ring of integers $R$ in
        $\mathbb{Q}(\zeta_m)$ is either $\mathcal{G}$ or
        $\mathcal{E}$.

        \item $\mu$ satisfies $INT$ or $\Sigma INT$

    \end{enumerate}
\end{enumerate}

\begin{Cor}\label{Cor:FiveSolutions-G}
There are five triples $(a,b,d)$, corresponding to $\nu =
(\frac{1}{a}, \frac{1}{b}, \frac{2}{d})$, which satisfy these
assumptions:
$$(3,2,12), (4,2,8), (6,2,6), (3,3,6), (4,4,4)$$
\end{Cor}
\begin{proof}
Because of property $G$, every point in the $\pi^{-1}(t_j)$ has
the same ramification index. By the Proposition, they are integer
multiples of $a, b,$ and $1$ respectively.  The assumptions
require $\mu_{0,i} = r_0\cdot \nu_0$ and $\mu_{1,i} = r_1 \cdot
\nu_1$ to all be integral, and furthermore $\sum_i \mu_{\infty, i}
= \frac{d}{r_{\infty}} \nu_{\infty} = 2$. In order for the ring
$R$ to be $\mathcal{G}$ or $\mathcal{E}$, the lowest common
denominator of the $\nu_i$ must be one of $2,3,4,6$. Finally,
using the fact $\sum_j \nu_j = 1$, one enumerates the solutions,
which yields the above list of five.  It is useful to note these
are precisely the solutions of
$$ \frac{1}{a} + \frac{1}{b} + \frac{2}{d} = 1. $$\label{eq:classify}

Because the non-integral $\mu_i$ are attached to $S =
\pi^{-1}(\infty)$, one sees the $\mu_i$ are all equal weight and
hence certainly satisfy $\Sigma INT$, and indeed correspond to one
of the ancestral examples or its equal weight descendants.
\end{proof}

Now we want to classify, under above the assumptions, when the
subspace $\Spi$ is actually a subball quotient.  By Theorem
\ref{Thm:MyMain}, given any valid $\nu$, this amounts to a
dimension count. The fact that $|T| = 3$ makes this easy to check.

\begin{Cor}\label{Cor:Five-are-Spi-as-BallQuot}
The triples $(a,b,d)$ from Corollary \ref{Cor:FiveSolutions-G} all
define $\Spi$ that are codimension $1$ subball quotients of
$DM(d,\mu)$, where $\mu_i = \frac{2}{d}, \forall i$.
\end{Cor}
\begin{proof}
Because $|T| = 3, dim_{\mathbb{C}}(IH_1(\pone, l_T)) = 1$.  It is
easy to check that the image of $\pi^*$ is non-trivial: simply
pull-back the generator $I_{1,\infty}$ and observe the resulting
linear combination of basis elements in the lift is never the
identity.  Thus the image of $\pi^*$ is one-dimensional.

By Theorem \ref{Thm:MyMain}, the $\Spi$ are subball quotients when
the dimension count agrees.  Here $im(\pi^*)^{\perp}$ is
codimension $1$.  Note that $|S| = d$. Then by Proposition
\ref{Prop:T-ram-dimcount} the codimension of $\Spi$ is determined
by the number of ramification conditions imposed by the
$T$-ramification class of $\pi$: precisely, $(d-3) - (2d-2 -
\frac{d}{a}(a-1) - \frac{d}{b}(b-1)) = (d-3) + 2 - d(\frac{1}{a} +
\frac{1}{b}) = d(1 - \frac{1}{a} - \frac{1}{b}) +2 - 3$, which,
using relation \ref{eq:classify} above, is simply $d(\frac{2}{d})
- 1 = 1$, as desired.  Note that one could also do a direct
dimension count from the explicit form of $\pi$ in this case.
\end{proof}

\begin{Rmk}\label{Rmk:ModuliRES}
The case (3,2,12) corresponds to the moduli space of rational
elliptic surfaces.  See Corollary
\ref{Cor:RationalEllipticSurfaces}.
\end{Rmk}

\subsubsection{Some moduli spaces of inhomogeneous forms and ball quotients}

Throughout let $A(u,v)$ and $B(u,v)$ be polynomials of degree
$d_1$ and $d_2$, respectively, with $d_1 < d_2$.  The definitions
can be extended to any number of polynomials, but we will use only
two.

\begin{Def}\label{Def:PseudoDisc}
Given the data $(A,B)$ as above, let $a, b$ and $N$ be positive
integers such that $d_1a = N = d_2 b$.  The choice of $a$ and $b$
determines a morphism $\Delta: \mathbb{W}\mathbb{P}^n(d_1,d_2)
\rightarrow \mathbb{P}^N$, given by $[A,B] \mapsto [A^a + B^b]$.
We call such a map a {\em pseudo-discriminant}.
\end{Def}
\begin{Rmk}
Note the map is well-defined.  Indeed, it is clear that the map
$(A,B) \mapsto A^a + B^b$ is $\mathbb{C}^*$-equivariant, where
$\mathbb{C}^*$ acts as multiplication by $(\lambda^{d_1},
\lambda^{d_2})$ and $\lambda^N$ respectively.
\end{Rmk}

\begin{Rmk}
One can interpret this map in the language of the GKZ theory of
resultants and discriminants for toric varieties \cite{GKZ}. There
it appears as an ``$A$-discriminant'', with $A$ an appropriately
chosen set of homogeneous polynomials, before quotienting out by
an associated group of toric automorphisms.
\end{Rmk}

The question we ask is essentially the following elementary (but
in many instances surprisingly rich) one.

\vspace{.1in}

\noindent {\bf Question:}  Given a degree $N$ polynomial in two
variables, when, and in how many ways, can it be written as the
sum of an $a^{th}$ power of a degree $d_1$ polynomial and a
$b^{th}$ power of a degree $d_2$ polynomial?

\vspace{.1in}

The ``when" is the image of $\Delta$ and the ``in how many ways"
is the degree of $\Delta$.  To be more precise, we are interested
in the number of solutions for a generic point in the image of
$\Delta$, not a complete analysis of the number of solutions for
any given degree $N$ polynomial.

Let $\zeta_m$ represent a primitive $m^{th}$ root of unity.  It is
clear that $\Delta([\zeta_a^j A, \zeta_b^k B]) = [A^a + B^b] =
\Delta([A,B])$. This is an obvious obstruction to the generic
injectivity of $\Delta$.

\begin{Prop} \label{Prop:Delta-atleast-gcd}
For a given $a$ and $b$, $\Delta$ is generically at least
$gcd(a,b)$-to-one.  In particular, for $\Delta$ to be generically
injective, it is necessary that $gcd(a,b) = 1$.
\end{Prop}
\begin{proof}
Because $\Delta$ is a map from weighted projective space
$\mathbb{W}\mathbb{P}^n(d_1,d_2)$ to projective space
$\mathbb{P}^N$, one must check which pairs $(\zeta_a^j A,
\zeta_b^k B)$ are equivalent under the weighted $\mathbb{C}^*$
action.

Clearly it suffices to check when there exists a complex number
$\lambda$ such that simultaneously $\lambda^{d_1} = \zeta_a^j$ and
$\lambda^{d_2} = \zeta_b^k$.  In particular, $\lambda$ must be an
$N^{th}$ root of unity. The two conditions are equivalent to
asking for solutions to the following system of congruences:
\begin{eqnarray*}
d_1 x \equiv d_1 j \pmod{N} &  \longleftrightarrow & x \equiv j \pmod{a} \\
d_2 x \equiv d_2 k \pmod{N} &  \longleftrightarrow & x \equiv k
\pmod{b}
\end{eqnarray*}
The Chinese Remainder Theorem implies there is a solution for all
$j$ and $k$ precisely when $a$ and $b$ are relatively prime.  More
generally, it implies that for any given $j$ there are
$N/gcd(a,b)$ values of $k$ which lie in the same
$\mathbb{C}^*$-orbit, so there are at least $N/(N/(gcd(a,b))) =
gcd(a,b)$ distinct points mapped to the same point by $\Delta$.
\end{proof}

What follows is a sufficient condition for the degree of the
pseudo-discriminant to be precisely $gcd(a,b)$.  Under this
circumstance, the sole obstruction to injectivity is the one
above, i.e., whether rescaling $(A,B)$ by relevant roots of unity
produces points in the same weighted $\mathbb{C}^*$-orbit.

\begin{Prop}\label{Prop:Delta-is-gcd}
Assume $d_2 > d_1 + 1$ and $b=2$.  Then $\Delta(A_1, B_1) =
\Delta(A_2, B_2)$ implies $A_1^a = A_2^a$ (equivalently, $B_1^b =
B_2^b$), that is, $A_1 = \zeta_a^j A_2$ and $B_1 = \zeta_b^k B_2$.
Furthermore, the degree of $\Delta$ is $gcd(a,b)$.
\end{Prop}
\begin{proof}
The argument we give is inspired by \cite[p. 17]{Vakil}. Consider
the space of polynomial quadruples $(A_1, B_1, A_2, B_2)$ subject
to the constraint that $A_1^a + B_1^b = A_2^a + B_2^b$. Remove the
subset $A_1^a = A_2^a$ (equivalently, $B_1^b = B_2^b$). What
remains are the solutions to $\Delta(A_1,B_1) = \Delta(A_2,B_2)$
other than $A_1 = \zeta_a^j A_2$ and $B_1 = \zeta_b^k B_2$.  Call
this set $Q_{\Delta}$. We claim $Q_{\Delta}$ is empty.

We argue by contradiction.  Assume it is not empty. Then it has a
dimension.  The dimension cannot be any less than the dimension of
the space of polynomials $(A_1, B_1)$, which is $(d_1 +1)+ (d_2 +
1) = d_1 + d_2 +2$.  So $dim_{\mathbb{C}}(Q_{\Delta}) \geq d_1 +
d_2 + 2$.

But there is another way to count the dimension of $Q_{\Delta}$.
Rewrite the defining constraint as $A_1^a - A_2^a = B_2^b -
B_1^b$. Because $b = 2$, the right hand side of the equation
factors as $(B_2-B_1)(B_2+B_1)$.  Specifying the pair $(A_1, A_2)$
determines $(B_2-B_1)(B_2+B_1)$.  Because $A_1^a \neq A_2^a$,
$A_1^a - A_2^a$ has $N$ roots (counting multiplicity).  By
assigning these roots to each of $(B_2 - B_1)$ and $(B_2 + B_1)$,
these factors are completely determined up to relative scaling and
the finite ambiguity in assigning the roots. Thus the dimension of
the set of solutions $Q_{\Delta}$ is the dimension of $(A_1, A_2)$
plus one to account for the relative scaling. That is,
$dim_{\mathbb{C}}(Q_{\Delta})  = 2(d_1+1) + 1 = 2d_1 + 3$.

One concludes $2d_1 + 3 \geq d_1 + d_2 + 2$, hence $d_1 + 1 \geq
d_2$, or equivalently $d_1 + 1 > d_2$.  But this contradicts the
assumption of our theorem that $d_2 > d_1 + 1$. Thus $Q_{\Delta}$
must be empty.
\end{proof}
\begin{Rmk}
Although the $b=2$ condition can be relaxed, the $d_2 > d_1 + 1$
is necessary.  As an example, when $d_1 = 2, d_2 = 3, N = 6$,
$\Delta$ is generically a $40$-to-$1$ map \cite{Elk}.
\end{Rmk}

\begin{Prop}
The pseudo-discriminant $\Delta$ is
$SL_2(\mathbb{C})$-equivariant.  It descends to a map of GIT
quotients.
\end{Prop}
\begin{proof}
The $SL_2(\mathbb{C})$-equivariance is immediate, because it acts
through the standard representation on $V \cong \mathbb{C}^2$ in
each case: the domain is $\mathbb{W}\mathbb{P}^n(Sym^{d_1}(V)
\oplus Sym^{d_2}(V))$ and the range is $\mathbb{P}(Sym^N(V))$.
\end{proof}

One should think of the image of this map as lying inside the
moduli space of $N$ unordered points.  The domain and range both
offer potentially different compactifications for the open set. In
particular, for $N = 12$ or $N=8$ the compactification of the
image is a Baily-Borel compactification for the Eisenstein or
Gaussian ancestral examples respectively.  Thus there is an
alternate compactification to the GIT compactification for certain
weighted projective space quotients.

Observe this gives alternate description of the $\Spi$.

\begin{Thm}\label{Thm:PseudoDisc-Hypersurface}
The classification of $\Spi$ satisfying property $G$ in Corollary
\ref{Cor:Five-are-Spi-as-BallQuot} is identical to the
classification of pseudo-discriminants with image a hypersurface
of codimension $1$.
\end{Thm}
\begin{proof}
This is simply a dimension count.  The condition that the image of
$\Delta$ be a hypersurface is the statement that $(d_1 + 1) +
(d_2+1) = N$, where $N = ad_1 = bd_2$.  Divide by $N$ to get
$\frac{1}{a} + \frac{1}{b} + \frac{2}{N} = 1$, which is the same
constraint as the one we discovered in the classification of $\pi$
with property $G$.
\end{proof}

\begin{Cor}\label{Cor:PairedBinaryForms}
The moduli spaces of inhomogeneous binary forms of bidegree
$(a,b)$, for $(a,b)$ taken from the list in the above theorem, are
branched covers of subball quotients of the corresponding $\DM$.
For one of the cases, $(3,2,12)$, $\Delta$ is an embedding on a
suitable open subset, and for two others, $(6,2,6)$ and $(4,2,8)$,
it is generically $2$-to-$1$.
\end{Cor}
\begin{Rmk}
This result parallels the statement that the ancestral examples,
thought of as moduli spaces of binary forms of degree $8$ and
$12$, are ball quotients.
\end{Rmk}

\begin{Cor}\label{Cor:RationalEllipticSurfaces}
 The moduli space of rational elliptic surfaces is a
ball quotient, in particular it is a hyperball quotient of the
Eisenstein ancestral example.
\end{Cor}
\begin{proof}
The example $(3,2,12)$ is the GIT moduli space of rational
elliptic surfaces presented as rational Weierstrass fibrations.
This GIT description of the moduli space was first discovered by
Miranda \cite{Mir}, following Mumford.
\end{proof}

\end{document}